%% file: paper_iso_div_arxiv.tex
\documentclass[11pt]{article}
\usepackage{amsmath}
\usepackage{amssymb}
\usepackage{graphicx}
\usepackage{color}
\usepackage{amsthm}
\usepackage{amsfonts}
\usepackage{amscd}
\usepackage{mathrsfs}
\usepackage{bm}
\usepackage{enumerate}
\usepackage{algorithmic, algorithm}

\newcommand{\bx}{\mathbf{x}}
\newcommand{\bn}{\mathbf{n}}
\newcommand{\by}{\mathbf{y}}

\newcommand{\bs}{\mathbf{s}}
\newcommand{\p}{\partial}
\newcommand{\al}{\alpha}

\newcommand{\diver}{\text{div}}

\newcommand{\hk}{R_t(\bx, \by)}
\newcommand{\rhk}{\bar{R}_t(\bx, \by)}
\newcommand{\rrhk}{\bar{\bar{R}}_t(\bx, \by)}

\newcommand{\hkxpj}{R_t(\bx, {\bf x}_j)}

\newcommand{\rhkxpj}{\bar{R}_t(\bx, {\bf x}_j)}

\newcommand{\bff}{{\bf f}}

\newcommand{\bfu}{{\bf u}}
\newcommand{\bfp}{{\bf x}}
\newcommand{\bfs}{{\bf s}}

\newcommand{\bz}{\mathbf{z}}
\newcommand{\bq}{\mathbf{q}}

\newcommand{\M}{{\mathcal M}}
\newcommand{\invt}{\frac{1}{t}}

\newcommand{\cM}{{\mathcal M}}

\newcommand{\bV}{\mathbf{V}}
\newcommand{\bA}{\mathbf{A}}

\newcommand{\mathd}{\mathrm{d}}

 \newtheorem{theorem}{\textbf{Theorem}}[section]
 
 \newtheorem{remark}{\textbf{Remark}}[section]
 \newtheorem{proposition}{\textbf{Proposition}}[section]
 \newtheorem{assumption}{\textbf{Assumption}}[section]
 \newtheorem{definition}{\textbf{Definition}}[section]

\newcommand{\R}{\mathbb{R}}

\numberwithin{equation}{section}

\makeatother

\begin{document}

\title{A convergent point integral method
for isotropic elliptic equations on point cloud \thanks{This research was supported by NSFC Grant 11201257 and 11371220.}}

\author{
Zhen Li%
\thanks{Yau Mathematical Sciences Center, Tsinghua University, Beijing, China,
100084. \textit{Email: zli12@mails.tsinghua.edu.cn.}%
} \and Zuoqiang Shi%
\thanks{Yau Mathematical Sciences Center, Tsinghua University, Beijing, China,
100084. \textit{Email: zqshi@math.tsinghua.edu.cn.}%
} }

 \maketitle


\begin{abstract}
In this paper, we propose a numerical method to solve isotropic elliptic equations on point cloud by
 generalizing the point integral method. The idea of the point integral method is to approximate the differential operators by integral operators
and discretize the corresponding integral equation on point cloud. The key step is to get the integral approximation.
In this paper, with proper kernel function, we get an integral approximation for the elliptic operators with isotropic coefficients. 
Moreover, the integral approximation has been proved to keep the coercivity of the original 
elliptic operator. The convergence of the point integral method is also proved.    
\end{abstract}




\section{Introduction}

Nowadays, data plays more and more important roles in science and engineering. 
In many problems, data is usually represented as a collection of
points embedding in a high dimensional Euclidean space. 
Processing and analysis of the point cloud data is essential in many applications, such as machine learning
\cite{belkin2003led, Coifman05geometricdiffusions} and image processing \cite{Peyre09,LDMM}.

In many applications, the point cloud data lies in a manifold whose 
dimension is much lower than the ambient Euclidean space. The low dimensionality is an important feature we could exploit to 
analyze the data.    
One example is the low dimensional manifold model (LDMM) in image processing \cite{LDMM}.
In this model, the original image is cutted to many overlap patches. The collection of all patches consists of a point cloud in 
Euclidean space. It is found that for many natural images, the patch set usually samples a low dimensional manifold which is called patch manifold. 
The dimension of the patch manifold is used as a regularization to processing the image. Based on differential geometry and variational method, 
this model is reduced to solve Laplace equation on patch set. The key point in LDMM is to solve this Laplace equation accurately and efficiently. 

Beside the data analysis, solving PDEs on manifold also appears in many physical problems, such as
material science \cite{CFP97,EE08}, fluid flow \cite{GT09,JL04},
biology and biophysics \cite{BEM11,ES10,NMWI11,WD08}.
To solve PDEs on manifold, many methods have been developed, especially on 2D surfaces, including 
surface finite element method \cite{DE-Acta}, level set method \cite{Bertalmio,XZ03}, grid based particle method \cite{LZ09,LZ11} 
and closest point method \cite{RM08,MR09}. However, these methods need extra information besides the point cloud, 
for instance, meshes, level set function and closest point function. These information is not easy to obtain from point cloud when the dimension of the 
manifold is high. 

Recently, Liang et al. proposed to discretize the differential operators on point cloud by local least square approximations of the manifold \cite{Liang13}.
Their method can achieve high
order accuracy and enjoy more flexibility since no mesh is needed. In principle, it can be applied to
manifolds with arbitrary dimensions and co-dimensions with or without boundary. However, if the dimension of the manifold is high, 
this method may not be stable since high order polynomial is used to fit the data. Later, Lai et al. proposed local mesh method to approximate the 
differential operators on point cloud \cite{Lai13}. The main idea is to construct mesh locally around each point by using K nearest neighbors. The local mesh is easier 
to construct than global mesh. Based on the local mesh, it is easy to discretize differential operators and compute integrals. However, when the dimension of the 
manifold is high, even local mesh is not easy to construct. 

The original point integral method for Laplace equation is closely related with the graph Laplacian \cite{Chung, pcdlp2009}.
Graph Laplacian has been widely used in many problems. 
It is observed in~\cite{BelkinN05, Lafon04diffusion, Hein:2005:GMW, Singer06} that
the graph Laplacian with the Gaussian weights well approximates the Laplace-Beltrami operator when the 
vertices of the graph are assumed to sample the underlying manifold.
When there is no boundary, Belkin and Niyogi~\cite{CLEM_08} showed the spectra of the graph Laplacian with 
Gaussian weights converges to that of Laplace-Beltrami operator. 
Recently, Singer and Wu~\cite{Singer13} showed the spectral convergence of the graph Laplacian in 
the presence of the Neumann boundary. 

Inspired by the graph Laplacian and the nonlocal diffusion, 
we developed the point integral method for Poisson equation on point cloud \cite{LSS,SS-neumann,SS-dirichlet}.
\begin{equation*}
      -\Delta_\mathcal{M} u(\bx)=f(\bx),\quad \bx\in \mathcal{M},
\end{equation*}
where $\Delta_\M=\diver(\nabla)$ is the Laplace-Beltrami operator in $\M$.

We assume that $\cM\in C^{\infty}$ is a compact $k$-dimensional manifold isometrically embedded in $\mathbb{R}^d$ with the standard Euclidean metric and $k\le d$.
If $\cM$ has boundary, the boundary, $\p\cM$ is also a $C^\infty$ smooth manifold.

Let $\Phi: \Omega\subset \mathbb{R}^k\rightarrow \cM\subset\mathbb{R}^d$ be a local parametrization of $\cM$ and $\theta\in \Omega$.
For any differentiable function $f:\cM\rightarrow \mathbb{R}$,
 define the gradient on the manifold
\begin{align}
  \label{eq:diff-M}
  \nabla f(\Phi(\theta))&=\sum_{i,j=1}^m g^{ij}(\theta)\frac{\p \Phi}{\p\theta_i}(\theta)\frac{\p f(\Phi(\theta))}{\p\theta_j}(\theta),
\end{align}
and for vector field $F:\M\rightarrow T_\bx\M$ on $\M$, where $T_\bx\M$ is the tangent space of $\M$ at $\bx\in \M$, the divergence is defined as
\begin{align}
\label{eq:diver}
\diver (F)&= \frac{1}{\sqrt{\det G}}\sum_{k=1}^d\sum_{i,j=1}^m\frac{\p}{\p \theta_i}\left(\sqrt{\det G}g^{ij}F^k(\Phi(\theta))\frac{\p \Phi^k}{\p\theta_j}\right)
\end{align}
where $(g^{ij})_{i,j=1,\cdots,k}=G^{-1}$, $\det G$ is the determinant of matrix $G$ and $G(\theta)=(g_{ij})_{i,j=1,\cdots,k}$ is the first fundamental form which is defined by
\begin{eqnarray}
  \label{eq:remainn}
  g_{ij}(\theta)=\sum_{k=1}^d\frac{\p \Phi_k}{\p\theta_i}(\theta)\frac{\p \Phi_k}{\p\theta_j}(\theta),\quad i,j=1,\cdots,m.
\end{eqnarray}
and $(F^1(\bx),\cdots,F^d(\bx))^t$ is the representation of $F$ in the embedding coordinates.

 The main idea of the point integral method is to
approximate the Poisson equation by the following integral equation:
\begin{equation*}
  -\int_\M  \Delta_\M u(\by)\rhk d\mu_\by\approx \invt\int_{\M} \hk(u(\bx) - u(\by))\mathd\mu_\by-2\int_{\p\M}\rhk \frac{\p u}{\p\bn}(\by)\mathd \tau_\by,
\end{equation*}
where $\bn$
is the out normal of $\M$, $R_t(\bx,\by)$ and $\bar{R}_t(\bx,\by)$ are kernel functions given as follows
\begin{equation}
\label{eq:kernel}
R_t(\bx, \by) = C_tR\left(\frac{|\bx -\by|^2}{4t}\right),\quad
\bar{R}_t(\bx, \by) = C_t\bar{R}\left(\frac{|\bx -\by|^2}{4t}\right)
\end{equation}
where $C_t = \frac{1}{(4\pi t)^{k/2}}$ is the normalizing factor.
 $R\in C^2(\mathbb{R}^+) $ be a positive function which is integrable over $[0,+\infty)$. And
\begin{equation*}
  \bar{R}(r)=\int_r^{+\infty}R(s)\mathd s.
\end{equation*}
There is not any derivatives in the integral equation. It is easy to
be discretized from point clouds using some quadrature rule.
In \cite{SS-neumann, SS-dirichlet}, we proved the convergence of the point integral method for Poisson equation with Neumann and Dirichlet boundary condition.

In the point integral method, we only need the point cloud to discretize the differential operator.
 This gives PIM great flexibility to fit the requirements in 
variety of applications.  
However, one limitation of the point integral method
is that it only applies on Laplace-Beltrami operator. In many problems, 
we need to discretize other differential operators besides Laplace-Beltrami operator. 
In this paper, we generalize the point integral method to isotropic elliptic operators. Isotropic elliptic operators are also widely used in many problems. 
One example is the nonlocal total variation minimization on 
point cloud, in which we need to solve an optimization problem,
\begin{align*}
  \min_{u} \|\nabla u\|_{L^1(\M)},\quad \text{subject to:} \quad \Psi(u)=b
\end{align*}
where $\nabla$ is the gradient in $\M$, $\Psi$ is the measurement operator related with the application, $b$ is the observation and 
\begin{align*}
  \|\nabla u\|_{L^1(\M)}=\int_\M |\nabla u (\bx)|\mathd \bx
\end{align*}
Using standard variational approach, the solution of above optimization problem can be given by solving a nonlinear elliptic equation, 
\begin{align*}
  -\diver\left(\frac{\nabla u(\bx)}{|\nabla u(\bx)|}\right)=f(\bx).
\end{align*}
where $f(\bx)$ is a known function.  
Apparently, this equation can be solved by solving a sequence of isotropic elliptic equation iteratively. 

In this paper, we consider to solve elliptic equations with isotropic coefficients on manifold $\cM$,
\begin{align}										
\label{eq:iso-intro}
      - \diver(p^2(\bx) \nabla u(\bx))  = f(\bx),\quad \bx\in \cM
\end{align}
The coeffcients $p(\bx)$ and source term $f(\bx)$ are
known smooth functions of spatial variables, i.e.
\begin{eqnarray*}
  p\in C^1(\cM), \; f \in C^1(\cM).
\end{eqnarray*}
The elliptic condition makes that there exist generic constants $c_0, c_1>0$ such that for any $\bx\in \M$,
\begin{eqnarray*}
  0<c_0\le p(\bx)\le c_1<\infty,
\end{eqnarray*}
The key observation in this paper is the integral approximation of isotropic elliptic operators given as following
\begin{align}
\label{eq:integral-iso}
 -\int_\M  \diver(p^2(\by)\nabla u(\by))\frac{\rhk}{p(\by)}d\mu_\by\approx &\invt\int_{\M} \hk(u(\bx) - u(\by))p(\by)\mathd\mu_\by\\
&-2\int_{\p\M} \frac{\p u}{\p\bn}(\by)\rhk p(\by)\mathd \tau_\by,\nonumber
\end{align}
where the kernel functions $R_t$ and $\bar{R}_t$ are same as those in \eqref{eq:kernel}. The main advantage of this integral approximation is that 
there is no differential operator inside. Using this approximation, we transfer the numerical differential to numerical integral which is much easier 
to compute on point cloud. 
Based on this integral approximation, we are able to develop the point integral method to isotropic elliptic equations. 

Similar integral approximation is also widely used in nonlocal diffusion and peridynamic model \cite{Du-SIAM,book-nonlocal,DGLZ13,DJTZ13,ZD10}.
The integral approximation is easy to implement on point cloud, since it has no derivatives inside.
Moreover, the point integral method also has very good theoretical property. It is proved that the coercivity of the original elliptic operator is
partially preserved and this partial coercivity implies the convergence of the point integral method.  


The rest of the paper is organized as following.
In Section 2, we introduce the point integral method for isotropic elliptic operator 
with Neumann and Dirichlet boundary condition. The convergence analysis is given in Section 3. 
Several numerical examples are presented in Section 4.
The conclusion remarks are made in Section 5.

\section{Point Integral Method for Isotropic Elliptic Equations}

In this section, we introduce a numerical method for isotropic elliptic equation on point cloud based on
the integral approximation \eqref{eq:integral-iso}. 

To simplify the notation, we introduce an integral operator,
\begin{equation}
  \label{eq:Lt}
L_tu(\bx)=\invt\int_{\M} \hk(u(\bx) - u(\by))p(\by)\mathd\mu_\by
\end{equation}
where $R_t$ is the kernel function given in \eqref{eq:kernel}.


\subsection{Neumann Boundary}
First, we consider the Neumann problem,
\begin{align}
  \label{eq:neumann}
  \left\{\begin{array}{cc}
-\diver(p^2(\bx)\nabla u(\bx))=f(\bx),& \bx\in \M,\\
\frac{\p u}{\p \bn}(\bx)=g(\bx),& \bx\in \p\M.
\end{array}\right.
\end{align}
Using the integral approximation \eqref{eq:integral-iso}, 
the solution of the Neumann problem \eqref{eq:neumann} can be obtained approximately by solving an integral equation
\begin{equation}
  \label{eq:integral-neumann}
  L_t u(\bx)-2\int_{\p\M} g(\by)\rhk p(\by)\mathd \tau_\by
=\int_\M  f(\by)\frac{\rhk}{p(\by)}d\mu_\by,
\end{equation}
with $t\ll 1$.

The eigenvalue problem is also solved by a generalized eigenvalue problem
\begin{equation}
  \label{eq:neumann-eigen}
  L_t u(\bx)=\lambda\int_\M  u(\by)\frac{\rhk}{p(\by)}d\mu_\by,
\end{equation}


\subsection{Dirichlet Boundary}
The Dirichlet problem is more involved in point integral method, since the normal derivative, $\frac{\p u}{\p \bn}$ is not known. 
\begin{equation}
  \left\{\begin{array}{rl}
     -\diver(p^2(\bx)\nabla u(\bx)) = f(\bx),& \bx\in \M,\\
      u(\bx) = g(\bx),&  \bx\in \p \M.
\end{array}\right.
\label{eq:dirichlet}
\end{equation}

Here, we use the same idea as that in \cite{LSS} to deal with the Dirichlet boundary. 

\subsubsection{Robin Approximation}
The simplest way is using Robin boundary to approximate the Dirichlet boundary. 
 More specifically, we consider the 
following Robin problem
\begin{equation}
  \left\{\begin{array}{rl}
     -\diver(p^2(\bx)\nabla u(\bx)) = f(\bx),& \bx\in \M,\\
      u(\bx)+\beta\frac{\p u}{\p \bn}(\bx) = g(\bx),&  \bx\in \p \M.
\end{array}\right.
\label{eq:robin}
\end{equation}
where $0<\beta\ll 1$ is a small parameter. It is easy to show that as $\beta \rightarrow 0$, the solution of the Robin problem, \eqref{eq:robin}, 
converges to the solution of the Dirichlet problem, \eqref{eq:dirichlet}. 

For the Robin problem, the integral approximation \eqref{eq:integral-iso} is applicable to give an integral equation,
\begin{align}
\label{eqn:integral_robin}
L_tu(\bx)
-&\frac{2}{\beta}\int_{\p \M} (g(\by)-u(\by)) \bar{R}_t(\bx,\by) p(\by)\mathd \tau_\by= \int_{\M}  \frac{f(\by)}{p(\by)}\bar{R}_t(\bx,\by) \mathd \mu_\by.
\end{align}
Similarly, we also get an approximation of the eigenvalue problem,
\begin{eqnarray}
\label{eqn:eigen_integral_robin}
L_tu(\bx)+\frac{2}{\beta}\int_{\p \M} u(\by) \bar{R}_t(\bx,\by)p(\by) \mathd \tau_\by
= \lambda \int_{\M} u(\by) \frac{\bar{R}_t(\bx,\by)}{p(\by)} \mathd \mu_\by
\end{eqnarray}

\subsubsection{Iterative Solver based on Augmented Lagrangian Multiplier}

In the Robin approximation, the parameter $\beta$ has to be small to get good approximation, while the linear system becomes ill-conditioned.
To alleviate this difficulty, we could use an iterative method based on the Augmented Lagrange method (ALM).

It is well known that the Dirichlet problem can be reformulated to be following
constrained optimization problem:
\begin{eqnarray}
\label{opt-dirichelet}
&&  \min_{v\in H^1(\mathcal{M})} \frac{1}{2}\int_\mathcal{M} p^2(\by) |\nabla v(\by)|^2 \mathd \mu_\by+\int_\mathcal{M} f(\by) v(\by) \mathd \mu_\by,\\
&& \mbox{subject to:} \; \;\;
v(\bx)|_{\p\mathcal{M}}=g(\bx), \nonumber
\end{eqnarray}

Applying the ALM method to the problem \eqref{opt-dirichelet}, 
we get an iterative method, in each step, an unconstrained optimization problem is solved,
\begin{eqnarray}
\label{alg-2}
\min_v \quad&&\frac{1}{2}\int_\mathcal{M} p^2(\by)|\nabla v(\by)|^2 \mathd \mu_\by+\int_\mathcal{M} f(\by) \cdot v(\by) \mathd \mu_\by\nonumber\\
&&+\int_{\p\mathcal{M}}w^k(\by)\cdot(g(\by)-v(\by))p^2(\by)\mathd \tau_\by
+\frac{1}{2\beta}\int_{\p\mathcal{M}}(g(\by)-v(\by))^2p^2(\by)\mathd \tau_\by.
\end{eqnarray}
Using the variational method, one can show that the solution to \eqref{alg-2}
is exactly the solution to the following Robin problem:
\begin{eqnarray}
\label{eq-robin-ALM}
  \left\{\begin{array}{cl}
      \diver(p^2(\bx) \nabla v(\bx))=f(\bx),& \bx\in \mathcal{M},\\
      v(\bx)+\beta \frac{\p v}{\p \bn}(\bx)=g(\bx)+\beta w^k(\bx),&  \bx\in \p \mathcal{M}.
\end{array}\right.
\end{eqnarray}
This Robin problem is solved by the integral equation. Notice that, the parameter $\beta$ is not necessarily small. Usually, we set $\beta=1$. Thus, the linear system is not
ill-conditioned. 
\begin{algorithm*}[t]
\floatname{algorithm}{Procedure}
\caption{ALM for Dirichlet Problem}
\label{alg-diri-approx}
\begin{algorithmic}[1]
\STATE $k=0$, $w^0=0$.
\REPEAT  
\STATE Solving the following integral equation to get $v^k$,
\begin{eqnarray*}
L_tv^k(\by)&-&\frac{2}{\beta}\int_{\p \mathcal{M}} (g(\by)-v^k(\by)+\beta w^k(\by)) \bar{R}_t(\bx,\by)p(\by) \mathd \tau_\by\
=\int_{\mathcal{M}} \frac{f(\by)}{p(\by)} \bar{R}_t(\bx,\by) \mathd \mu_\by.
\end{eqnarray*}
\STATE $w^{k+1}=w^k+\frac{1}{\beta}(g-(v^k|_{\p\M}))$, $k = k+1$
\UNTIL {$\|g-(v^{k-1}|_{\p \mathcal{M}})\| = 0$}
\STATE $u=v^k$
\end{algorithmic}
\end{algorithm*}


\subsection{Discretization}
The main advantage of the integral equations is that they are easy to discretize over the point cloud since there is not any derivatives inside. 

Assume we are given
a set of sample points $P$ sampling the submanifold $\M$ and a subset $S\subset P$ sampling the boundary of $\M$.
List the points in $P$ respectively $S$ in a fixed
order $P=(\bx_1, \cdots, \bx_n)$ where $\bfp_i \in \R^d, 1\leq i\leq n$, respectively $S=(\bfs_1, \cdots, \bfs_{n_b})$ where $\bfs_i \in P$.
In addition, assume we are also given
two vectors $\bV = (V_1, \cdots, V_n)^t$ where $V_i$ is an volume weight of $\bx_i$ in $\M$, and
$\bA= (A_1, \cdots, A_{n_b})^t$ where $A_i$ is an area weight of $\bfs_i$ in $\p \M$. In this point cloud data $(P,S, \bV,\bA)$, the integral
equation \eqref{eq:integral-neumann} can be disretized as
\begin{equation}
\invt\sum_{\bx_j \in P} R_t(\bx_i,\bx_j)(u_i - u_j)p_jV_j  = 2\sum_{\bs_j \in S} \bar{R}_t(\bx_i,\bs_j) b_j p_j A_j + 
\sum_{\bx_j \in P} \bar{R}_t(\bx_i,\bx_j) f_jV_j/p_j.
\label{eqn:dis}
\end{equation}
where $p_j=p(\bx_j),\;  f_j=f(\bx_j),\;b_l=b(\bs_l),\; j=1,\cdots,|P|,\; l=1,\cdots,|S|$.

The other integral equations and corresponding eigenvalue problems can be discretized consequently.

\begin{remark}
The integral approximation \eqref{eq:integral-iso} also holds if the parameter $t$ depends on $\bx$, i.e.
\begin{align}
\label{eq:integral-adaptive}
 -\int_\M  \mbox{\rm div}(p^2(\by)\nabla u(\by))\frac{\rhk}{p(\by)}d\mu_\by\approx &\frac{1}{t(\bx)}\int_{\M} R\left(\frac{|\bx-\by|^2}{4t(\bx)}\right)
(u(\bx) - u(\by))p(\by)\mathd\mu_\by\\
&-2\int_{\p\M} \frac{\p u}{\p\bn}(\by)\bar{R}\left(\frac{|\bx-\by|^2}{4t(\bx)}\right) p(\by)\mathd \tau_\by.\nonumber
\end{align}
Based on above approximation, in the computation, we can choose $t$ adaptive to the distribution of the points.
\end{remark}

\section{Convergence Analysis}

In this section, we analyze the convergence of the point integral method for isotropic elliptic equation. 
To make the theoretical analysis concise, we only consider the homogeneous Neumann boundary
conditions,
\begin{align}
  \label{eq:neumann-homo}
  \left\{\begin{array}{cc}
-\diver(p^2(\bx)\nabla u(\bx))=f(\bx),& \bx\in \M,\\
\frac{\p u}{\p \bn}(\bx)=0,& \bx\in \p\M.
\end{array}\right.
\end{align}
The corresponding numerical scheme is 
\begin{equation}
\invt\sum_{\bx_j \in P} R_t(\bx_i,\bx_j)(u_i - u_j)p_jV_j  = \sum_{\bx_j \in P} \bar{R}_t(\bx_i,\bx_j) f_jV_j/p_j.
\label{eqn:dis-homo}
\end{equation}
 The analysis can be easily generalized to
the non-homogeneous boundary conditions. The convergence of Dirichlet problem can be proved also following the similar procedure as that in 
\cite{SS-dirichlet}.

\subsection{Main Result}
\label{sec:assume}
We will prove that the solution given by the point integral method converges to the exact solution as the point cloud $(P,\mathbf{V})$
converges to the manifold $\M$. Before giving the result of the convergence, we need to clarify the meaning of the convergence of the point cloud $(P,\mathbf{V})$
to the manifold $\M$. 

First, we introduce an index to measure the distance between the point cloud $(P,\mathbf{V})$ and the manifold $\M$, which is called {\it integral accuracy index},
denoted as $h(P,\mathbf{V},\M)$. 
\begin{definition}[Integral Accuracy Index]
  \label{def:h}
For the point cloud $(P,\mathbf{V})$ which samples the manifold $\M$, the integral accuracy index $h(P,\mathbf{V},\M)$ is defined as
\begin{equation*}
h(P,\mathbf{V},\M)=\sup_{f\in C^1(\M)}\frac{\left|\int_\M f(\by) \mathd\mu_\by - \sum_{\bfp_i\in P} f(\bfp_i)V_i\right|}{|\text{\rm supp}(f)|\|f\|_{C^1(\M)}}.
\end{equation*}
where $\|f\|_{C^1(\M)} = \|f\|_\infty +\|\nabla f\|_\infty$ and
$|\text{\rm supp}(f)|$ is the volume of the support of $f$.
\end{definition}

Using the definition of integrable index, we say that the point cloud $(P,\mathbf{V})$ converges to the manifold $\M$ if 
$h(P,\mathbf{V},\M)\rightarrow 0$. In the convergence analysis, we consider the case that $h(P,\mathbf{V},\M)$ is small enough.
\begin{remark}
  In some sense, $h(P,\mathbf{V},\M)$ is a measure of the density of the point cloud. If the point cloud is uniformly distributed on the manifold, 
from central limit theorem, $h(P,\mathbf{V},\M)\sim O(1/\sqrt{|P|})$ where $|P|$ is the number of point in $P$.
\end{remark}
\begin{remark}
  To consider the non-homogeneous Neumann boundary condition or Dirichlet boundary condition, we have to also assume that $h(S,\mathbf{A},\p\M)\rightarrow 0$,
where $S$ is the point set sample the boundary $\p\M$ and $\mathbf{A}$ is the corresponding volume weight on the boundary $\p\M$.
\end{remark}

To get the convergence, we also need some assumptions on the regularity of the submanifold $\M$ and
the integral kernel function $R$.
\begin{assumption}
\label{assumptions}
\begin{itemize}
\item \rm Smoothness of the manifold: $\M, \p\M$ are both compact and $C^\infty$ smooth $k$-dimensional submanifolds isometrically embedded in a Euclidean space $\mathbb{R}^d$.
\item \rm Ellipticity: there exist generic constants $c_0,c_1>0$, such that $c_0\le p(\bx)\le c_1$ and $p(\bx)\in C^1(\M)$.
\item \rm Assumptions on the kernel function $R(r)$:
\begin{itemize}
\item[\rm (a)] \rm Smoothness: $R\in C^2(\mathbb{R}^+)$;
\item[(b)] Nonnegativity: $R(r)\ge 0$ for any $r\ge 0$.
\item[(c)] Compact support:
$R(r) = 0$ for $\forall r >1$;
\item[(d)] Nondegeneracy:
 $\exists \delta_0>0$ so that $R(r)\ge\delta_0$ for $0\le r\le\frac{1}{2}$.
\end{itemize}
\end{itemize}
\end{assumption}
\begin{remark}
  The assumption on the kernel function is very mild. The compact support assumption can be relaxed to exponentially decay, like Gaussian kernel.
 In the nondegeneracy assumption, $1/2$ may be replaced by a positive number $\theta_0$ with $0<\theta_0<1$. Similar assumptions on the kernel function 
is also used in analysis the nonlocal diffusion problem \cite{DLZ13}.
\end{remark}

All the convergence analysis in this paper is based on above assumptions. 
In the statement of the theorems, above assumptions are omitted to make the statements more concise. 

The other issue we have to address is that how to compute the difference between the discrete solution and the analytic solution.
The solution of the discrete system \eqref{eqn:dis-homo} is a vector $\bf u$ defined on $P$ while
the solution of the problem~\eqref{eq:neumann-homo} is a function defined on $\M$. To make them
comparable, for any solution $\bfu = (u_1, \cdots, u_n)^t,\; n=|P|$ to the problem~\eqref{eqn:dis-homo},
we construct a function on $\M$
\begin{equation}
\label{eqn:interp_neumann}
I_{\bff}(\bfu) (\bx) = \frac{ \sum_{\bx_j \in P} \hkxpj u_j p_j V_j + t\sum_{\bfp_j \in P} 
\rhkxpj f_jV_j/p_j} {\sum_{\bfp_j \in P} \hkxpj p_jV_j }.
\end{equation}
It is easy to verify that $I_{\bff}(\bfu)$ interpolates $\bfu$ at the sample points $P$, i.e.,
$I_{\bff}(\bfu)(\bfp_j) = u_j$ for any $\bx_j\in P$.
The following theorem guarantees the convergence of the point integral method.
\begin{theorem}
Let
$u$ be the solution to problem~\eqref{eq:neumann-homo} with  $f\in C^1(\M)$ and let
the vector $\bfu$ be the solution to the problem~\eqref{eqn:dis-homo}. Then there exists constants $C$ and $T_0$
depend on $\M$ and $p(\bx)$,
such that for any $t, \frac{h(P,\mathbf{V},\M)}{\sqrt{t}}\le T_0$,
\begin{equation}
  \|u-I_{\bff}(\bfu)\|_{H^1(\M)} \leq C\left(t^{1/2} + \frac{h(P,\mathbf{V},\M)}{t^{3/2}}\right)\|f\|_{C^1(\M)}.
\end{equation}
where $h(P,\mathbf{V},\M)$ is the integral accuracy index. 
\label{thm:poisson_neumann}
\end{theorem}

\subsection{Proof of Convergence}
\label{sec:intermediate}
Roughly,
 the proof the convergence includes two parts: estimate of the truncation error $L_t(u-I_{\bff}(\bfu))$
and the stability of the integral operator $L_t$. Here $L_t$ is the integral operator in \eqref{eq:Lt}
, $u(\bx)$ is the solution
of the problem~\eqref{eq:neumann-homo} and $\bfu$ is the solution of
the problem~\eqref{eqn:dis-homo}.

This strategy is 
standard in numerical analysis. It is well known that consistency together with stability imply convergence.
On the other hand, the point integral method has some 
special structures both in truncation error and stability, which makes the analysis a little more involved. 

First, we have following theorem regarding the stability of the operator $L_t$. 
\begin{theorem}
Let $u(\bx)$ solves the integral equation
\begin{eqnarray*}
  L_t u = r(\bx)
\end{eqnarray*}
where $r\in H^1(\M)$ with $\int_\M r(\bx)p(\bx)\mathd\mu_\bx=0$.
Then, there exist constants $C>0, T_0>0$ independent on $t$, such that
\begin{eqnarray*}
  \|u\|_{H^1(\M)}\le C\left(\|r\|_{L^2(\M)}+t\|\nabla r\|_{L^2(\M)}\right)
\end{eqnarray*}
as long as $t\le T_0$.
\label{thm:regularity}
\end{theorem}

To use above stability result, we need $L_2$ estimate of $L_t(u-I_{\bff}(\bfu))$ and $\nabla L_t(u-I_{\bff}(\bfu))$.
In the analysis, we split the truncation error $L_t(u-I_{\bff}(\bfu))$ to two terms,
\begin{equation*}
  L_t(u-I_{\bff}(\bfu))=L_t(u-u_t))+L_t(u_t-I_{\bff}(\bfu))
\end{equation*}
where $u_t$ is the solution of the integral equation 
\begin{equation}
  \label{eq:integral-homo}
\invt\int_{\M} \hk(u(\bx) - u(\by))p(\by)\mathd\mu_\by=\int_\M f(\by)\frac{\bar{R}_t(\bx,\by)}{p(\by)}\mathd \mu_\by.
\end{equation}

For the second term, we have following estimate.
\begin{theorem}
Let $u_t(\bx)$ be the solution of the problem~~\eqref{eq:integral-homo} and $\bfu$ be the solution of
the problem~\eqref{eqn:dis-homo}. If
$f\in C^1(\M)$ ,
 then there exists constants
$C, T_0$ depending only on $\M$ and the coefficient $ p(\bx)$, so that
\begin{eqnarray}
\|L_{t} \left(I_{\bff}\bfu - u_{t}\right)\|_{L^2(\M)} &\leq& \frac{Ch(P,\mathbf{V},\M)}{t^{3/2}}\|f\|_{C^1(\M)},\\
\label{eqn:dis_error_l2}
\|\nabla L_{t} \left(I_{\bff}\bfu - u_{t}\right)\|_{L^2(\M)} &\leq& \frac{Ch(P,\mathbf{V},\M)}{t^{2}}\|f\|_{C^1(\M)}.
\label{eqn:dis_error_dl2}
\end{eqnarray}
as long as $t\le T_0$ and $\frac{h(P,\mathbf{V},\M)}{\sqrt{t}}\le T_0$, $h(P,\mathbf{V},\M)$ is the integral difference index in 
Definition \ref{def:h}.
\label{thm:dis_error}
\end{theorem}

The error term $L_t(u-u_t))$ is a little more complicated. It has two parts, one is the interior term and the other is the boundary term. 
We need to estimate these two terms separately to get better estimation of the convergence rate.  
\begin{theorem}
 Let $u(\bx)$ be the solution
of the problem~\eqref{eq:neumann-homo} and $u_t(\bx)$ be the solution of the corresponding
integral equation \eqref{eq:integral-homo}. Let
\begin{align}
\label{eq:error_boundary}
&I_{bd} =\sum_{j=1}^d \int_{\p\M}n^j(\by)(\bx-\by)\cdot\nabla(\nabla^ju(\by))\rhk p(\by)\mathd \tau_\by,
\end{align}
and
\begin{align}
L_t (u- u_t)=I_{in}+I_{bd}\nonumber.
\end{align}
where $\bn(\by)=(n^1(\by),\cdots,n^d(\by))$ is the out normal vector of $\p\M$ at $\by$, $\nabla^j$ is the $j$th component of gradient $\nabla$.

If $u\in H^3(\M)$, then there exists constants
$C, T_0$ depending only on $\M$ and $p(\bx)$, so that,
\begin{eqnarray}
\label{eqn:integral_error_int}
\left\|I_{in}\right\|_{L^2(\M)}\le Ct^{1/2}\|u\|_{H^3(\mathcal{M})},\quad
\left\|\nabla I_{in}\right\|_{L^2(\M)} \leq C\|u\|_{H^3(\mathcal{M})},
\end{eqnarray}
as long as $t\le T_0$.
\label{thm:integral_error}
\end{theorem}

Using the definition of the boundary term $I_{bd}$, \eqref{eq:error_boundary},
it is easy to check that 
\begin{eqnarray*}
\left\|I_{bd}\right\|_{L^2(\M)}= O(t^{1/4}),\quad
\left\|\nabla I_{bd}\right\|_{L^2(\M)} = O(t^{-1/2}),
\end{eqnarray*}
Based on this estimation, Theorem \ref{thm:regularity} and Theorem \ref{thm:integral_error} give that 
\begin{equation*}
  \|u-u_t\|_{H^1(\M)}=O(t^{1/4}).
\end{equation*}
This proves the convergence, however the convergence rate is relatively low. This low rate comes from the boundary term. From interior term only, 
the rate is $\sqrt{t}$. Notice that the boundary term has a specific integral formula given in \eqref{eq:error_boundary}. 
Using this formula, we know that the boundary term concentrates in a small layer adjacent to the boundary whose width is of the order of $\sqrt{t}$ and vanish in the 
interior region.
Utilizing this special structure, we could get better convergence rate with the help of a stability estimate specifically for the boundary term, which is 
given in Theorem \ref{thm:regularity_boundary}.

\begin{theorem}
Let $u(\bx)$ solves the integral equation
\begin{eqnarray*}
  L_t u = \int_{\p\M}\mathbf{b}(\by)\cdot(\bx-\by)\rhk p(\by)\mathd \tau_\by-\bar{b}
\end{eqnarray*}
where $|\M|_p=\int_\M p(\bx)\mathd \mu_\bx$ and 
\begin{eqnarray*}
  \bar{b}=\frac{1}{|\M|_p}\int_\M \left(\int_{\p\M}\mathbf{b}(\by)\cdot(\bx-\by)\rhk p(\by)\mathd \tau_\by\right)\mathd\bx.
\end{eqnarray*}
Then, there exist constant $C>0, T_0>0$ independent on $t$, such that
\begin{eqnarray*}
  \|u\|_{H^1(\M)}\le C\sqrt{t}\;\|\mathbf{b}\|_{H^1(\M)}.
\end{eqnarray*}
as long as $t\le T_0$.
\label{thm:regularity_boundary}
\end{theorem}

Based on above four theorems, it is easy to prove Theorem~\ref{thm:poisson_neumann}. Using Theorem \ref{thm:dis_error} and Theorem \ref{thm:regularity},
we get $$\|u_t-I_{\bff}(\bfu)\|_{H^1(\M)}=O\left(\frac{h(P,\mathbf{V},\M)}{t^{3/2}}\right).$$ 
Applying Theorem \ref{thm:regularity} to the interior term in Theorem \ref{thm:integral_error} and Theorem \ref{thm:regularity_boundary} to the boundary term
respectively, we have $$\|u-u_t\|_{H^1(\M)}=O\left(t^{1/2}\right).$$ Putting above two inequality together, Theorem~\ref{thm:poisson_neumann} is proved. 

Next, we prove Theorem 
  \ref{thm:regularity}, \ref{thm:dis_error}, \ref{thm:integral_error} and \ref{thm:regularity_boundary} respectively.

\subsection{Proof of Theorem \ref{thm:integral_error}}
\label{sec:integral_error}

\input{integral_eqn}

\subsection{Proof of Theorem \ref{thm:dis_error}}

\input{discretization}

\subsection{Proof of Theorem \ref{thm:regularity}}

In order to prove Theorem \ref{thm:regularity}, we need two theorems, \ref{thm:elliptic_v} and \ref{thm:elliptic_L_t}. 
The proof of these two theorems can be obtained by making minor revision of 
 the proof of Theorem 4.4 and 4.5 in \cite{SS-neumann}, the details of the proof are put in the supplementary material. 
\begin{theorem} For any function $u\in L^2(\mathcal{M})$,
there exists a constant $C>0$ independent on $t$ and $u$, such that
  \begin{eqnarray}
    \int_\M\int_{\mathcal{M}}R\left(\frac{|\bx-\by|^2}{4t}\right) (u(\bx)-u(\by))^2p(\bx)p(\by)\mathd\mu_\bx\mathd \mu_\by \ge C\int_\mathcal{M} |\nabla v(\bx)|^2\nonumber
p(\bx)\mathd \mu_\bx,
  \end{eqnarray}
where 
\begin{eqnarray}
v(\bx)=\frac{C_t}{w_t(\bx)}\int_{\mathcal{M}}R\left(\frac{|\bx-\by|^2}{4t}\right) u(\by)p(\by)\mathd \mu_\by,\nonumber
\label{eqn:smooth_v}
\end{eqnarray}
and $w_t(\bx) = C_t\int_{\mathcal{M}}R\left(\frac{|\bx-\by|^2}{4t}\right)p(\by)\mathd \mu_\by$.
\label{thm:elliptic_v}
\end{theorem}

\begin{theorem}
Assume both $\M$ and $\p\M$ are $C^\infty$. There exists a constant $C>0$ independent on $t$
so that for any function $u\in L_2(\M)$ with $\int_\M u(\bx)p(\bx)\mathd \mu_\bx = 0$ and for any sufficient small $t$
\begin{eqnarray}
\int_\M\int_{\mathcal{M}}R\left(\frac{|\bx-\by|^2}{4t}\right) (u(\bx)-u(\by))^2p(\bx)p(\by)\mathd\mu_\bx\mathd \mu_\by \geq C\|u\|_{L_2(\M)}^2.\nonumber
\end{eqnarray}
\label{thm:elliptic_L_t}
\end{theorem}

Using above two theorems, Theorem \ref{thm:regularity} becomes an easy corollary.
\begin{proof}{\it of Theorem \ref{thm:regularity}}

Using Theorem \ref{thm:elliptic_L_t}, we have
  \begin{eqnarray}
\label{eqn:stable_Lt_l2}
    \|u\|_{L^2(\M)}^2 &\le&C\int_\M u(\bx)r(\bx)p(\bx)\mathd\mu_\bx\le C\|u\|_{L^2(\M)}\|r\|_{L^2(\M)}.
  \end{eqnarray}
This inequality \eqref{eqn:stable_Lt_l2} implies that
\begin{eqnarray*}
  \|u\|_{L^2(\M)}\le C\|r\|_{L^2(\M)}.
\end{eqnarray*}
Now we turn to estimate $\|\nabla u\|_{L^2(\M)}$. 
Notice that we have the following expression for $u$,
\begin{eqnarray*}
u(\bx)=v(\bx)+\frac{t}{w_t(\bx)}\, r(\bx)
\end{eqnarray*}
where
\begin{eqnarray*}
  v(\bx)=\frac{1}{w_t(\bx)}\int_{\M}R_t(\bx,\by)u(\by)p(\by)\mathd\mu_\by,\quad w_t(\bx)=\int_{\M}R_t(\bx,\by)p(\by)\mathd\mu_\by.
\end{eqnarray*}
By Theorem \ref{thm:elliptic_v}, we have
\begin{eqnarray*}
\|\nabla u\|_{L^2(\M)}^2&\le&   2\|\nabla v\|_{L^2(\M)}^2+ 2t^2\left\|\nabla \left(\frac{r(\bx) - \bar{r}}{w_t(\bx)}\right)\right\|_{L^2(\M)}^2\\
&\le & C \int_\M u(\bx)L_tu(\bx)p(\bx)\mathd \mu_\bx + Ct\|r\|_{L^2(\M)}^2+Ct^2\|\nabla r\|_{L^2(\M)}^2\nonumber\\
&\le& C\|u\|_{L^2(\M)}\|r\|_{L^2(\M)}+ Ct\|r\|_{L^2(\M)}^2+Ct^2\|\nabla r\|_{L^2(\M)}^2\nonumber\\
&\le& C\|r\|_{L^2(\M)}^2+Ct^2\|\nabla r\|_{L^2(\M)}^2\nonumber\\
&\le & C\left(\|r\|_{L^2(\M)}+t\|\nabla r\|_{L^2(\M)}\right)^2.\nonumber
\end{eqnarray*}
The proof is completed.   
\end{proof}

\subsection{Proof of Theorem \ref{thm:regularity_boundary}}

\begin{proof}
First, we denote
\begin{align*}
r(\bx)&=\int_{\p\M}\mathbf{b}(\by)\cdot(\bx-\by)\rhk p(\by)\mathd \tau_\by,\\
\bar{r}&=\frac{1}{|\M|_p}\int_\M \left(\int_{\p\M}\mathbf{b}(\by)\cdot(\bx-\by)\rhk p(\by)\mathd \tau_\by\right)p(\bx)\mathd\bx.
\end{align*}
where $|\M|_p=\int_\M p(\by)\mathd \mu_\by$.

The key point of the proof is to show that
  \begin{eqnarray}
\label{eq:est_boundary_whole}
   \left|\int_\M u(\bx)\left(r(\bx)-\bar{r}\right)p(\bx)\mathd \mu_\bx\right|
	\le C\sqrt{t} \;\|\mathbf{b}\|_{H^1(\M)} \|u\|_{H^1(\M)}.
  \end{eqnarray}

First, notice that
$$|\bar{r}|\le 
C\sqrt{t}\;\|\mathbf{b}\|_{L^2(\p\M)}\le C\sqrt{t}\;\|\mathbf{b}\|_{H^1(\M)}.$$
Then it is sufficient to show that
  \begin{equation}
\label{eq:est_boundary}
    \left|\int_\M u(\bx)\left(\int_{\p\M}\mathbf{b}(\by)\cdot(\bx-\by)\bar{R}_t(\bx,\by)
p(\by)\mathd \tau_\by\right)p(\bx)\mathd \mu_\bx\right|\le C\sqrt{t} \;\|\mathbf{b}\|_{H^1(\M)} \|u\|_{H^1(\M)}.
  \end{equation}
Direct calculation gives that
\begin{eqnarray*}
  &&|2t\nabla\rrhk-(\bx-\by)\bar{R}_t(\bx,\by)|\le C|\bx-\by|^2\rhk,
\end{eqnarray*}
where $\rrhk=C_t\bar{\bar{R}}\left(\frac{\|\bx-\by\|^2}{4t}\right)$ and $\bar{\bar{R}}(r)=\int_{r}^{\infty}\bar{R}(s)\mathd s$.
This implies that
\begin{align}
\label{eq:est_boundary_1}
  &\left|\int_\M u(\bx)p(\bx)\int_{\p\M}\mathbf{b}(\by)\left((\bx-\by)\bar{R}_t(\bx,\by)+2t\nabla\rrhk\right)
p(\by)\mathd \tau_\by \mathd \mu_\bx\right|\\
\le & C\int_\M |u(\bx)p(\bx)|\int_{\p\M}|\mathbf{b}(\by)||\bx-\by|^2\rhk p(\by)\mathd \tau_\by\mathd \mu_\bx\nonumber\\
\le &Ct\|\mathbf{b}\|_{L^2(\p\M)} \left(\int_{\p\M}\left(\int_\M\rhk p(\bx)\mathd\mu_\bx\right)
\left(\int_\M |u(\bx)|^2\rhk p(\bx)\mathd \mu_\bx\right)p(\by)\mathd \tau_\by\right)^{1/2}\nonumber\\
\le & Ct\|\mathbf{b}\|_{H^1(\M)} \left(\int_{\M} |u(\bx)|^2p(\bx)
\left(\int_{\p\M} \rhk p(\by)\mathd \tau_\by\right)\mathd \mu_\bx\right)^{1/2}\nonumber\\
\le & Ct^{3/4}\|\mathbf{b}\|_{H^1(\M)}\|u\|_{L^2(\M)}.\nonumber
\end{align}
On the other hand, using the Gauss integral formula, we have
 \begin{eqnarray}
\label{eq:gauss_boundary}
&&    \int_\M u(\bx)p(\bx)\int_{\p\M}\mathbf{b}(\by)\cdot\nabla\rrhk p(\by)\mathd \tau_\by\mathd \mu_\bx\\
&=& \int_{\p\M} \int_{\M}u(\bx)p(\bx)T_\bx(\mathbf{b}(\by))\cdot\nabla\rrhk p(\by)\mathd \mu_\bx \mathd \tau_\by\nonumber\\
&=&\int_{\p\M} \int_{\p\M}\mathbf{n}(\bx)\cdot T_\bx(\mathbf{b}(\by))u(\bx)\rrhk p(\bx)p(\by)\mathd \tau_\bx\mathd \tau_\by\nonumber\\
&&-\int_{\p\M} \int_{\M}\text{div}_\bx[u(\bx)p(\bx)T_\bx(\mathbf{b}(\by))]\rrhk p(\by)\mathd \mu_\bx\mathd \tau_\by.\nonumber
  \end{eqnarray}
Here $T_\bx$ is the projection operator to the tangent space on $\bx$. To get the first equality, we use the 
fact that $\nabla\rrhk$ belongs to the tangent space on $\bx$, such that $\mathbf{b}(\by)\cdot\nabla\rrhk=T_\bx(\mathbf{b}(\by))\cdot\nabla\rrhk$
and $\bn(\bx)\cdot T_\bx(\mathbf{b}(\by))=\bn(\bx)\cdot \mathbf{b}(\by)$ where $\bn(\bx)$ is the out normal of $\p\M$ at $\bx\in \p\M$.

For the first term, we have
  \begin{align}
\label{eq:est_boundary_2}
&    \left|\int_{\p\M} \int_{\p\M}\mathbf{n}(\bx)\cdot T_\bx(\mathbf{b}(\by))u(\bx)\rrhk p(\bx)p(\by)\mathd \tau_\bx\mathd \tau_\by\right|\\
= &\left|\int_{\p\M} \int_{\p\M}\mathbf{n}(\bx)\cdot \mathbf{b}(\by)u(\bx)\rrhk p(\bx)p(\by)\mathd \tau_\bx\mathd \tau_\by\right|\nonumber\\
\le &C\|\mathbf{b}\|_{L^2(\p\M)}  \left(\int_{\p\M} \left(\int_{\p\M}|u(\bx)|\rrhk p(\bx)\mathd \tau_\bx\right)^2
p(\by)\mathd \tau_\by\right)^{1/2}\nonumber\\
\le &C\|\mathbf{b}\|_{H^1(\M)}  \left(\int_{\p\M} \left(\int_{\p\M}\rrhk p(\bx)\mathd \tau_\bx\right)
 \left(\int_{\p\M}|u(\bx)|^2\rrhk p(\bx)\mathd \tau_\bx\right)
p(\by)\mathd \tau_\by\right)^{1/2}\nonumber\\
\le &Ct^{-1/2}\; \|\mathbf{b}\|_{H^1(\M)} \|u\|_{L^2(\p\M)}\le Ct^{-1/2}\; \|\mathbf{b}\|_{H^1(\M)} \|u\|_{H^1(\M)}.\nonumber
  \end{align}
We can also bound the second term on the right hand side of \eqref{eq:gauss_boundary}. By using the assumption that $\M\in C^\infty$, we have
\begin{align*}
&|  \text{div}_\bx[u(\bx)p(\bx)T_\bx(\mathbf{b}(\by))]|\nonumber\\
\le& |\nabla u(\bx)||T_\bx(\mathbf{b}(\by))||p(\bx)|
+|u(\bx)||\text{div}_\bx[T_\bx(\mathbf{b}(\by))]||p(\bx)|+|\nabla p(\bx)||u(\bx)T_\bx(\mathbf{b}(\by))| \\
\le & C(|\nabla u(\bx)|+|u(\bx)|)|\mathbf{b}(\by)|
\end{align*}
where the constant $C$ depends on the curvature of the manifold $\M$.

Then, we have
  \begin{eqnarray}
\label{eq:est_boundary_3}
&&    \left|\int_{\p\M} \int_{\M}\text{div}_\bx[u(\bx)T_\bx(\mathbf{b}(\by))]\rrhk p(\bx)p(\by)\mathd \mu_\bx\mathd \tau_\by\right|\\
&\le &C   \int_{\p\M}\mathbf{b}(\by)p(\by) \int_{\M}(|\nabla u(\bx)|+|u(\bx)|)\rrhk p(\bx)\mathd \mu_\bx\mathd \tau_\by\nonumber\\
&\le &C \|\mathbf{b}\|_{L^2(\p\M)}  \left(\int_{\M} (|\nabla u(\bx)|^2+|u(\bx)|^2)p(\bx)
\left(\int_{\p\M}\rrhk p(\by)\mathd \tau_\by\right)\mathd \mu_\bx\right)^{1/2}\nonumber\\
\quad\quad&\le &C t^{-1/4}\;\|\mathbf{b}\|_{H^1(\M)}  \|u\|_{H^1(\M)}.\nonumber
  \end{eqnarray}
Then, the inequality \eqref{eq:est_boundary} is obtained from \eqref{eq:est_boundary_1},
\eqref{eq:gauss_boundary}, \eqref{eq:est_boundary_2} and \eqref{eq:est_boundary_3}.
Now, using Theorem \ref{thm:elliptic_L_t}, we have
  \begin{eqnarray}
    \|u\|_{L^2(\M)}^2\le C \int_\M u(\bx)L_tu(\bx)p(\bx)\mathd\mu_\bx
    \le C\sqrt{t}\;\|\mathbf{b}\|_{H^1(\M)}  \|u\|_{H^1(\M)}.
\label{eq:est_l2_boundary}
  \end{eqnarray}
Note $r(\bx)=\int_{\p\M}(\bx-\by)\cdot \mathbf{b}(\by)\bar{R}_t(\bx,\by)p(\by)\mathd \tau_\by$. Direct calculation gives us that
\begin{eqnarray*}
  \|r(\bx)\|_{L^2(\M)}&\le& Ct^{1/4}\|\mathbf{b}\|_{H^1(\M)} ,~\text{and}\\
  \|\nabla r(\bx)\|_{L^2(\M)}&\le& Ct^{-1/4}\|\mathbf{b}\|_{H^1(\M)} .
\end{eqnarray*}
The integral equation $L_t u=r-\bar{r}$ gives that
\begin{eqnarray*}
u(\bx)=v(\bx)+\frac{t}{w_t(\bx)}\,(r(\bx)-\bar{r})
\end{eqnarray*}
where
\begin{eqnarray*}
  v(\bx)=\frac{1}{w_t(\bx)}\int_{\M}R_t(\bx,\by)u(\by)p(\by)\mathd\mu_\by,\quad w_t(\bx)=\int_{\M}R_t(\bx,\by)p(\by)\mathd\mu_\by.
\end{eqnarray*}
By Theorem \ref{thm:elliptic_v}, we have
\begin{eqnarray}
\label{eq:est_dl2_boundary}
&&\|\nabla u\|_{L^2(\M)}^2 \\
&\le&   2\|\nabla v\|_{L^2(\M)}^2+ 2t^2\left\|\nabla \left(\frac{r(\bx)-\bar{r}}{w_t(\bx)}\right)\right\|_{L^2(\M)}^2\nonumber\\
&\le & C \int_\M u(\bx)L_tu(\bx)p(\bx)\mathd\mu_\bx + Ct\|r\|_{L^2(\M)}^2+Ct^2\|\nabla r\|_{L^2(\M)}^2\nonumber\\
&\le& C\sqrt{t}\;\|\mathbf{b}\|_{H^1(\M)}  \|u\|_{H^1(\M)}+ Ct\|r\|_{L^2(\M)}^2+Ct^2\|\nabla r\|_{L^2(\M)}^2\nonumber\\
&\le& C\|\mathbf{b}\|_{H^1(\M)} \left(\sqrt{t}\|u\|_{H^1(\M)}+ Ct^{3/2}\right).\nonumber
\end{eqnarray}
Using \eqref{eq:est_l2_boundary} and \eqref{eq:est_dl2_boundary}, we have
\begin{eqnarray*}
  \|u\|_{H^1(\M)}^2\le C\|\mathbf{b}\|_{H^1(\M)}  \left(\sqrt{t}\|u\|_{H^1(\M)}+ Ct^{3/2}\right), 
\end{eqnarray*}
which proves the theorem.
\end{proof}

\section{Numerical Experiments}

In this section, we show several numerical examples to demonstrate the performance of the point integral method for isotropic elliptic equations. 
This section is separated to two parts. In the first part, on some simple 2D surfaces, the convergence of the point integral method is verified. 
In the second part, we consider a nonlocal total variation minimization problem, in which some isotropic elliptic equations are solved on point cloud in 
high dimensional space.

\subsection{Examples on 2D Surfaces}
 
%

In this subsection, we consider the isotropic elliptic equation on 2D surfaces
 \begin{align}
\label{eq:iso-num}
-\diver (p^2(\bx)\nabla u(\bx)) = f(\bx), \quad  \bx \in \cM,
\end{align}
with Neumann and Dirichlet boundary conditions,
\begin{align*}
\frac{\p u}{\p\bn}(\bx) = b(\bx), \quad \mbox{or}\quad u(\bx) = b(\bx),\quad  \bx \in \p\cM
\end{align*}
To estimate the volume weight vector $\bV$ from the point sets $P$, 
a local mesh around each sample point is constructed, from which the weight of that point is computed. For details to estimate the volume weight, 
we refer to \cite{LSS}. The kernel function is chosen to be Gaussian function,
\begin{align*}
  R_t(\bx,\by)=\frac{1}{(4\pi t)^{k/2}}\exp\left(-\frac{\|\bx-\by\|^2}{4t}\right).
\end{align*}
The parameter $t$ is set as $t = \left(\frac{1}{|P|}\sum_{i=1}^{|P|}\rho(\bx_i)\right)^2$, where $\rho(\bx_i)$ is the radius of 10 nearest neighbors of $\bx_i$.

\begin{figure}[!ht]
\begin{center}
\begin{tabular}{cc}
\includegraphics[width=0.45\textwidth]{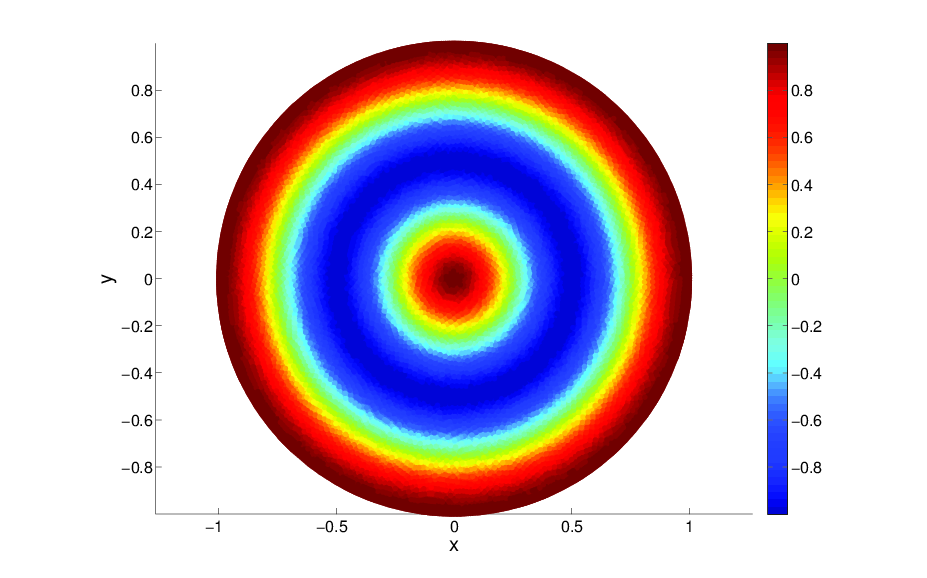}&
\includegraphics[width=0.45\textwidth]{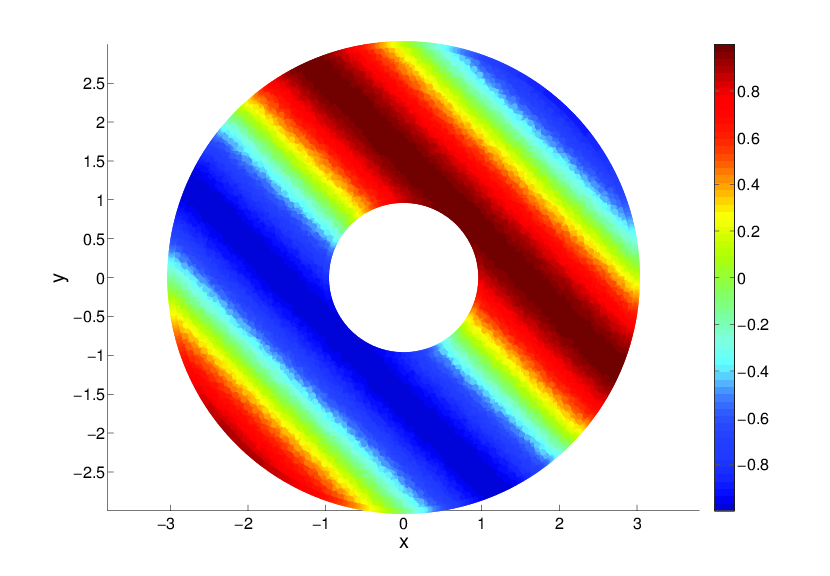}\\
(a)&(b)
\end{tabular}
\end{center}
\vspace{-4mm}
\caption{Ground truth: (a) $u_{gt} = \cos(2\pi r)$ in unit disk; (b) $u_{gt} = \sin(x+y)$ in the annulus}   \label{fig:other_u_gt}
\end{figure}

\paragraph{Example 1}
In the first example, the manifold $\M$ is an unit disk and an annulus in $\mathbb{R}^2$. The inner radius of the annulus is $1$ and outer radius is $3$.
The exact solution is set to be $u_{gt}(\bx) = \cos(2\pi \|\bx\|)$ in unit disk and $u_{gt} = \sin(x+y)$ in the annulus, see Figure \ref{fig:other_u_gt}.
The coefficient of the equation in \eqref{eq:iso-num} is
\begin{align}                                                                                   \label{coef_iso}
p = 1+\frac{1}{4}\|\bx\|^2,
\end{align}
both in the unit disk and annulus.
The Neumann boundary condition is enforced in unit disk and we consider the Dirichlet boundary condition in the annulus.

Table \ref{tbl:disk_Poisson} list the $l_2$ error of the point integral method as the number of points grows. 
This result clearly shows the convergence of the point integral method. The convergence rate in $l_2$ error is approximately $1/\sqrt{|P|}$. 
\begin{table}
\begin{center}
\begin{tabular}{| c| c | c | c | c |}
\hline
$|P|$   &	684		&	2610		&	10191		&	40296 \\
\hline
disk	&	0.364597	&	0.214960	&	0.111961	&	0.056028\\
\hline
annulus	&	0.036760	&	0.012227	&	0.005557	&	0.003542\\
\hline
\end{tabular}
\end{center}
\vspace{3mm}
\caption{$l_2$ error for $u_{gt} = \cos(2\pi r)$ in the unit disk and $\sin(x+y)$ in the annulus. \label{tbl:disk_Poisson}}	
\end{table}

The eigenvalue problem with homogeneous Neumann boundary condition is also solved in the annulus.
\begin{align*}
-\diver (p^2(\bx)\nabla u(\bx)) =& \lambda u(\bx), \quad \bx \in \cM \\
\frac{\p u}{\p\bn}(\bx) =& 0, \quad \quad\quad \bx \in \p\cM
\end{align*}
and the coefficient $p$ is given in \eqref{coef_iso}.

\begin{figure}[!ht]
\begin{center}
\includegraphics[width=0.6\textwidth]{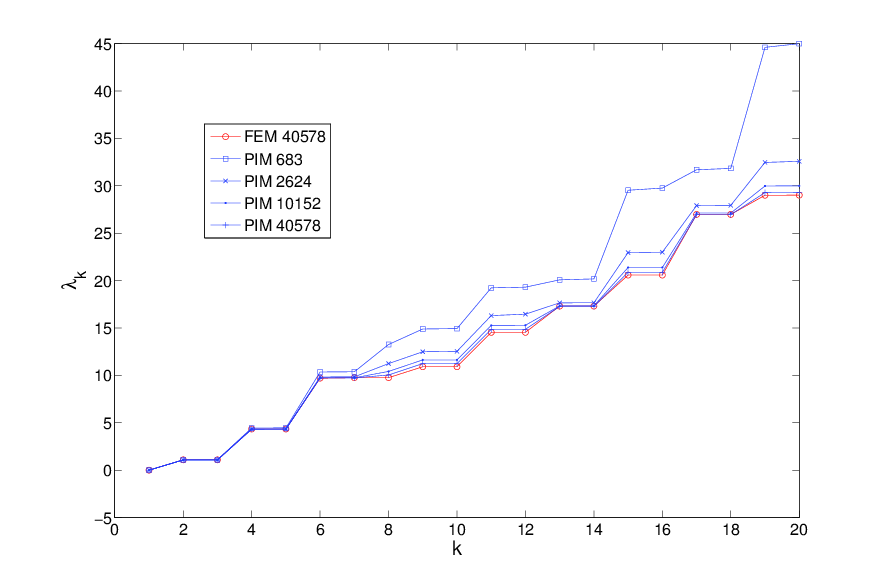}
\end{center}
\vspace{-4mm}
\caption{First 20 eigenvalues in the annulus with Neumann boundary condition with different point cloud.}			\label{fig:annulus_eigenvalues}
\end{figure}

The first 20 eigenvalues are plotted in Figure \ref{fig:annulus_eigenvalues}. The eigenvalues given by finite element method in the finest mesh 
is used as the true solution. Our result shows that the eigenvalue computed in the point integral method also converge.

\begin{figure}[!ht]
\begin{center}
\begin{tabular}{cc}
\includegraphics[width=0.6\textwidth]{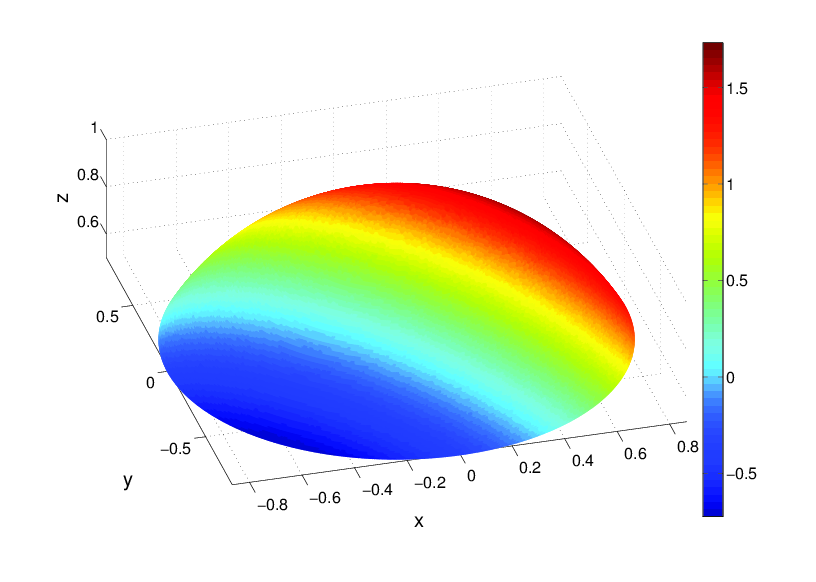}
\end{tabular}
\end{center}
\vspace{-4mm}
\caption{Ground truth: $u_{gt} = x+y+z$ on the cap.}				\label{fig:cap}
\end{figure}

\paragraph{Example 2}
Now, we solve equation \eqref{eq:iso-num} with Neumann condition and Dirichlet condition
on a curved surface in $\mathbb{R}^3$.
Let $\cM$ be a cap on the unit sphere, whose height is $1/2$ and the cap angle is $\pi/3$, as shown in Figure \ref{fig:cap}. 
The coefficient of the equation is also given in \eqref{coef_iso}. 

We set the ground truth to be $u_{gt} = x+y+z$, where $(x,y,z)$ is the coordinate in $\mathbb{R}^3$.

The $l_2$ errors of the point integral method are listed in Table \ref{tbl:cap_Poisson_x+y+z}.
The convergence rate for both boundary value problems are $1$.
\begin{table}[!ht]
\begin{center}
\begin{tabular}{| c| c | c | c | c |}
\hline
$|P|$	&	1199		&	4689		&	18540		&	73757\\
\hline
Neumann	&	0.036779	&	0.015355	&	0.007479	&	0.003189\\
\hline
Dirichlet&	0.007238	&	0.001921	&	0.001278	&	0.000750\\
\hline
\end{tabular}
\end{center}
\vspace{3mm}
\caption{$l_2$ error for $u_{gt} = x+y+z$ on cap.}								\label{tbl:cap_Poisson_x+y+z}
\end{table}

\begin{figure}[!ht]
\begin{center}
\includegraphics[width=0.6\textwidth]{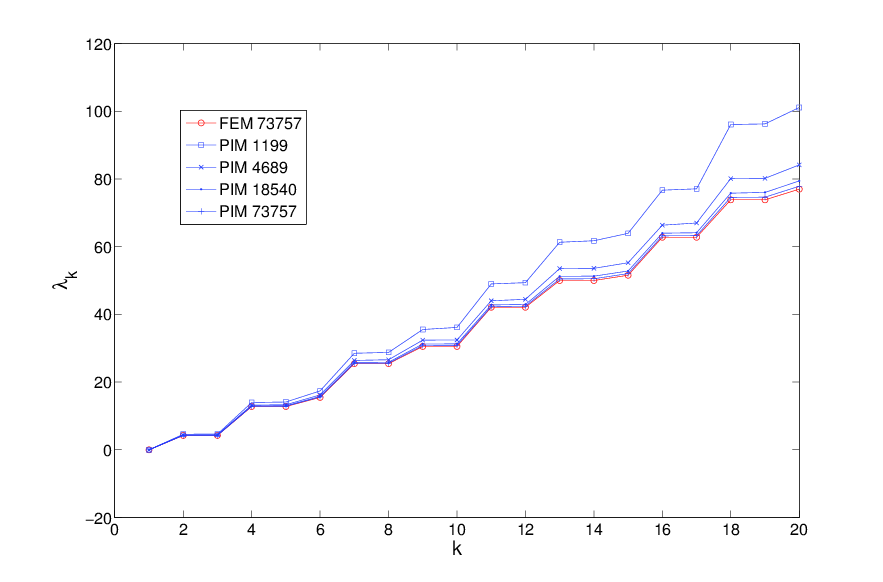}
\end{center}
\vspace{-4mm}
\caption{First 20 eigenvalues on the cap with Neumann boundary condition.}			\label{fig:cap_eigenvalues}
\end{figure}

The first 20 eigenvalues are also computed for homogeneous Neumann condition as shown in Figure \ref{fig:cap_eigenvalues}.
 As the number of points increases, the eigenvalues given by PIM converge to those computed by FEM, which suggests the convergence of the point integral method.

\begin{figure}[!ht]
\begin{center}
\begin{tabular}{cc}
\includegraphics[width=0.45\textwidth]{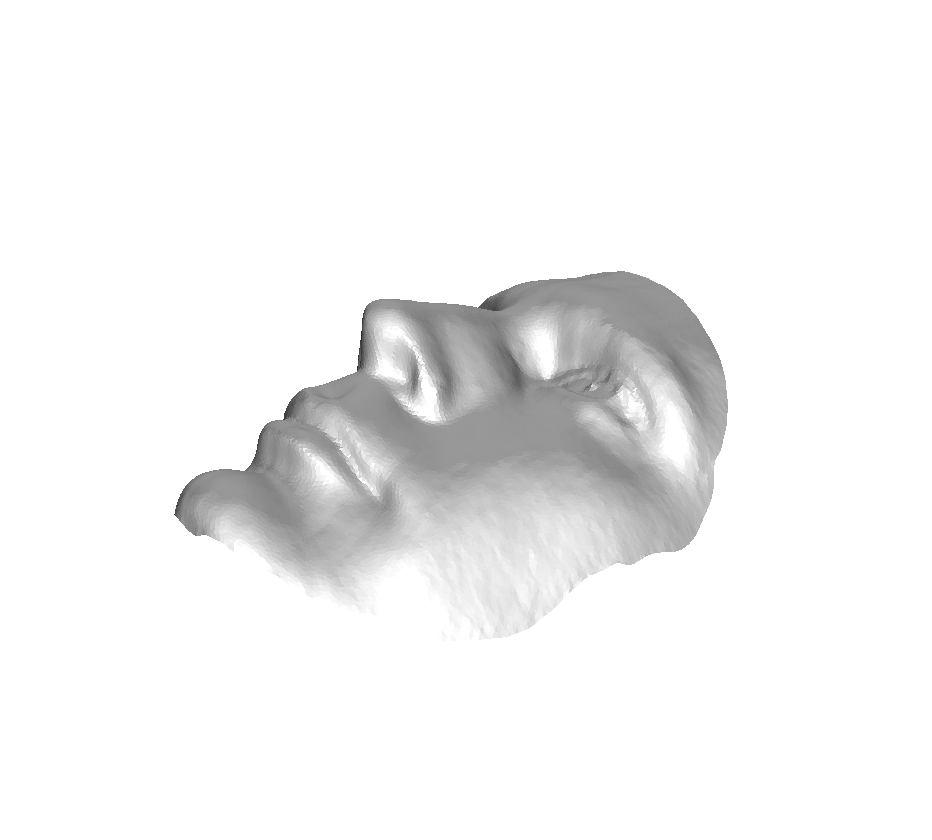}&
\includegraphics[width=0.45\textwidth]{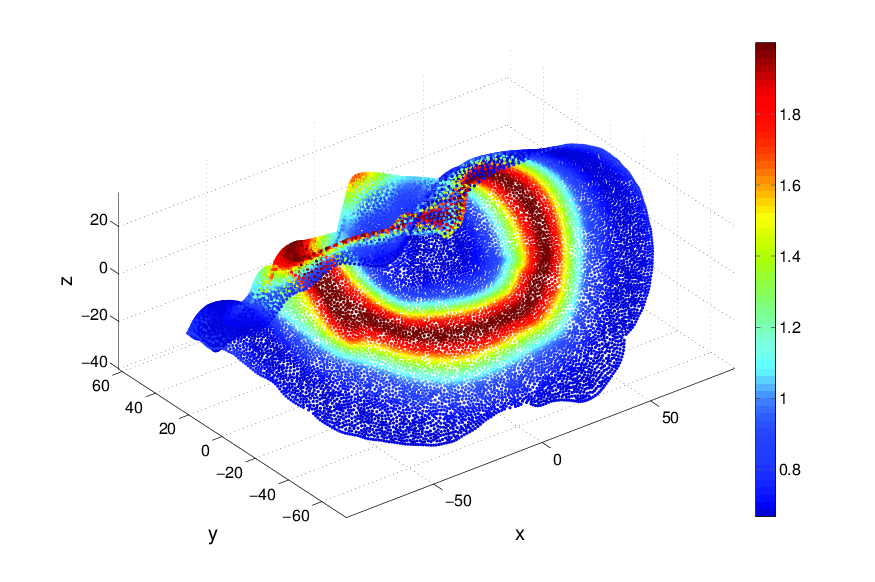}
\end{tabular}
\end{center}
\vspace{-4mm}
\caption{(a) Face of Alex; (b) Coefficient: restriction of $p^2$ on $\cM$}		\label{fig:Alex}
\end{figure}

\begin{figure}[!ht]
\begin{center}
\begin{tabular}{cc}
\includegraphics[width=0.45\textwidth]{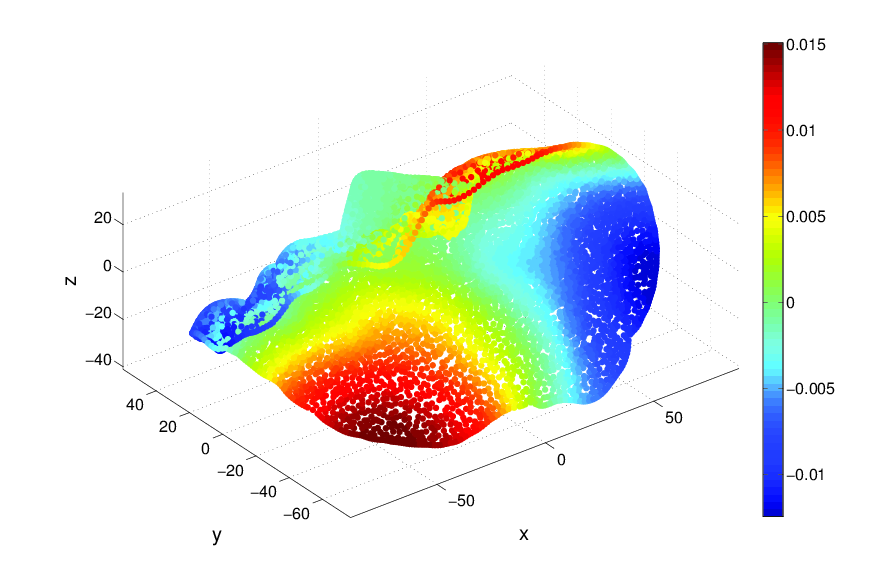}&
\includegraphics[width=0.45\textwidth]{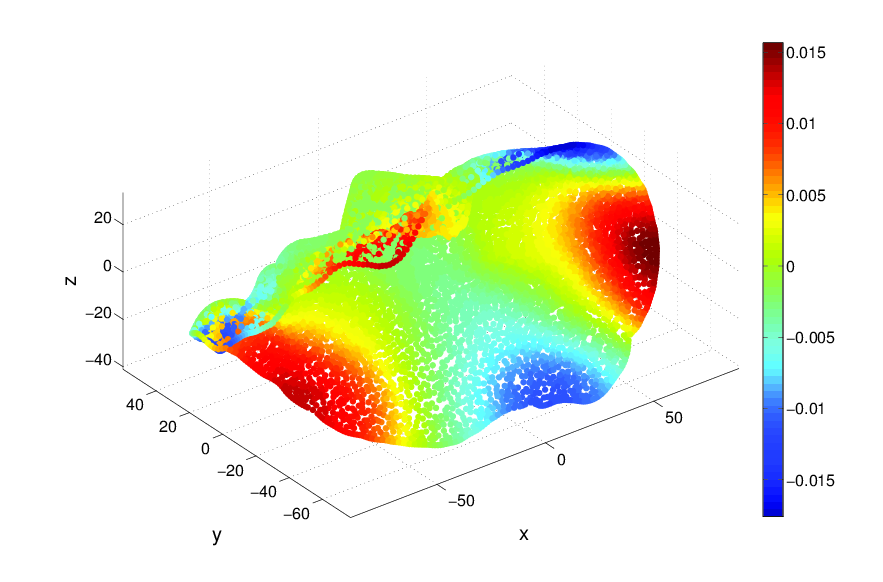}\\
\includegraphics[width=0.45\textwidth]{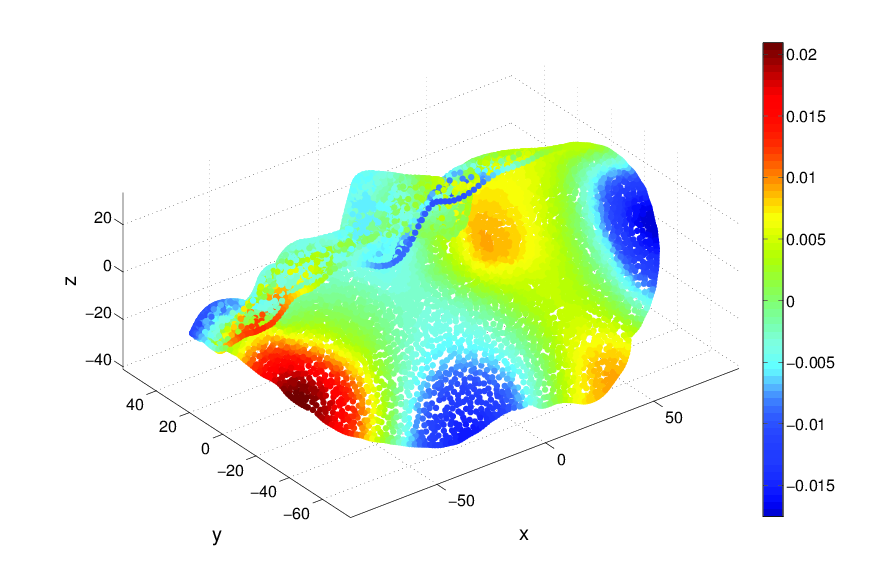}&
\includegraphics[width=0.45\textwidth]{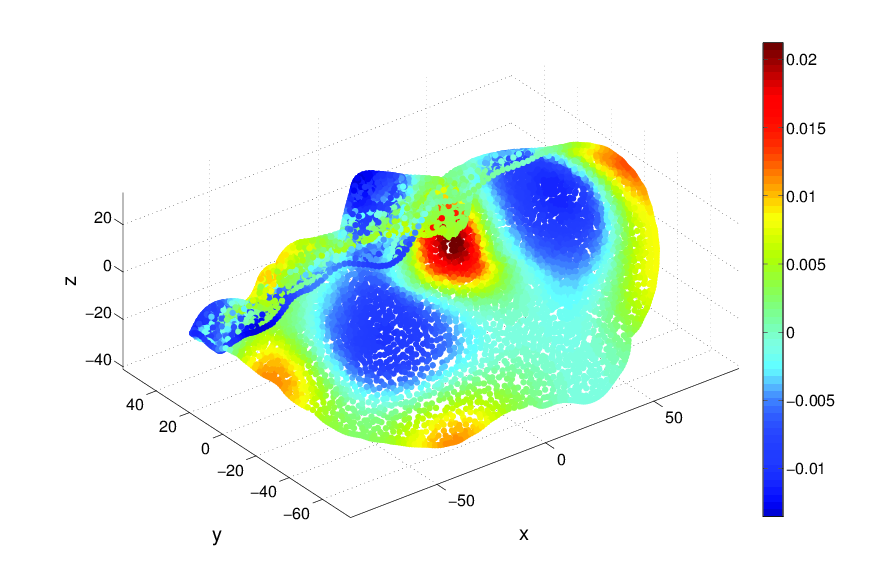}\\
\end{tabular}
\end{center}
\vspace{-4mm}
\caption{Eigenfunctions on `Alex' with homogenous Neumann boundary condition.}		\label{fig:Alex_eigenfunctions}
\end{figure}

\paragraph{Example 3}
In this example, we consider a more complex surface, a human face called "Alex". 
The surface is sampled by 10597 points (Figure \ref{fig:Alex}) and the analytic form of the surface is not known.
The coefficient of the equation in \eqref{eq:iso-num} is
\begin{align*}
p^2 = \frac{1}{\sin(r/10)/2+1}
\end{align*}
where $r = \sqrt{x^2+y^2+z^2}$. 
In this example, we solve the eigenvalue problem of the isotropic elliptic operator. 
Several eigenfunctions computed by the point integral method are shown in Figure \ref{fig:Alex_eigenfunctions}.

From the examples in 2D surfaces, we see that PIM solves isotropic elliptic equations with Neumann and Dirichlet boundary very well. Moreover, the 
convergence rate is higher than that obtained in the convergence analysis. The point integral method is applicable to point cloud in high dimensional space,
not only on the 2D surfaces. Next, we will show a high dimensional example.

\subsection{Nonlocal Total Variation Extension}

\begin{figure}[!ht]
\begin{center}
\begin{tabular}{cc}
\includegraphics[width=0.45\textwidth]{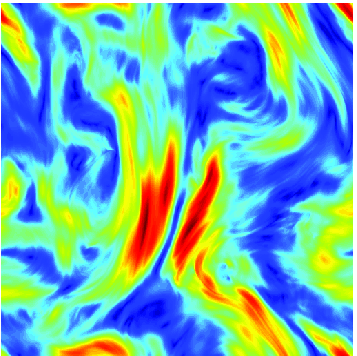}&
\includegraphics[width=0.45\textwidth]{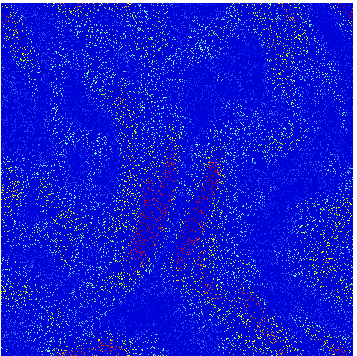}\\
(a)&(b)
\end{tabular}
\end{center}
\vspace{4mm}
\caption{(a): original data; (b): 10\% subsampled data.}   \label{fig:plasma}
\end{figure}

In this example, we consider an $L_1$ extension on point cloud. The point cloud is constructed by using the 
patches of a $512\times 512$ image, which is shown in  Figure \ref{fig:plasma}(a).
The original image is subsampled and only retain 10\% of the pixels at random. The subsampled image is shown in Figure \ref{fig:plasma}(b).
One classical problem in image processing is to recover the image from the subsampled image. Here, rather than 
give an image reconstruction method, we only use this
example to demonstrate the performance of the point integral method for isotropic elliptic equations.

In this example, the point cloud consists of the patches of the original image. 
For each pixel $x_i$ in the image $f$, we extract a patch around it of size $5\times 5$
which is denoted as $p_{x_i}(f)$, where $f$ is the original image. 
Totally, we get $512^2$ patches and each patch is $5\times 5$. The collection of all the patches give 
a point cloud in $\mathbb{R}^{25}$. Denote this point cloud as $P=\{p_{x_i}(f): i=1,\cdots, 512^2\}$.
The image is actually corresponding a function $u$ on the point cloud $P$ with $u(p_{x_i}(f))=f(x_i)$, $f(x_i)$ is the value
of image $f$ at pixel $x_i$. Corresponding to the subsampled image, the value of function $u$ is only known in the patches around 
the sampled pixels. The collection of all these patches is denoted as $S$. 

Recently, manifold model attracts many attentions in image processing \cite{LDMM}. In manifold model, the point cloud $P$ is assumed to be a sample of an 
underlying manifold, which is called patch manifold. The total variation is used as a regularization to reconstruct the image. 
The main idea is to minimize the total variation 
in the patch manifold, i.e.,
\begin{align}
  \label{eq:l1-ext}
  \min_u \|\nabla u\|_{L^1(\M)},\quad \mbox{subject to:}\quad u(\bx)=f(\bx),\quad \bx\in S. 
\end{align}
The variation approach tells us that the optimal solution of \eqref{eq:l1-ext} is given by solving following PDE,
\begin{align*}
  \diver\left(\frac{\nabla u(\bx)}{|\nabla u(\bx)|}\right)=0,
\end{align*}
with the Dirichlet type boundary condition
\begin{align*}
   u(\bx)=f(\bx),\quad \bx\in S.
\end{align*}
One natural method to solve above PDE is an iterative scheme,
\begin{align}
\label{eq:l1-ite}
  \diver\left(\frac{\nabla u^{n+1}(\bx)}{|\nabla u^{n}(\bx)|}\right)=0, \quad \quad \quad u^{n+1}(\bx)=f(\bx),\quad \bx\in S.
\end{align}
In each step, we need to solve an isotropic elliptic equation.

Here, the gradient is computed by using an integral approximation also.
\begin{align*}
  \nabla u(\bx)=\frac{1}{t \,\bar{w}_t(\bx)}\int_{\M}R_t(\bx,\by)(\bx-\by)(u(\bx)-u(\by))\mathd \mu_\by
\end{align*}
$\bar{w}_t(\bx)=\int_\M \bar{R}_t(\bx,\by)\mathd \mu_\by$.
In the computation, to avoid degenerate of the ellipticity, we regularize the coefficient by adding a small constant in the denominator, i.e., replace
$|\nabla u^{n}(\bx)|$ by $|\nabla u^{n}(\bx)|+\epsilon$ in \eqref{eq:l1-ite} with $\epsilon=10^{-3}$.
The point cloud is assumed to be uniformly distributed, so the volume weight is uniform. The kernel function is Gaussian function. In this example,
we use the integral approximation \eqref{eq:integral-adaptive} with adaptive $t(\bx_i) = \rho(\bx_i)^2$, where $\rho(\bx_i)$ is the radius of 20 nearest neighbors of $\bx_i$.


\begin{figure}[!ht]
\begin{center}
\begin{tabular}{cc}
\includegraphics[width=0.45\textwidth]{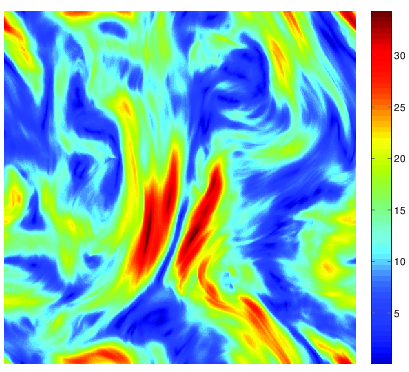}&
\includegraphics[width=0.45\textwidth]{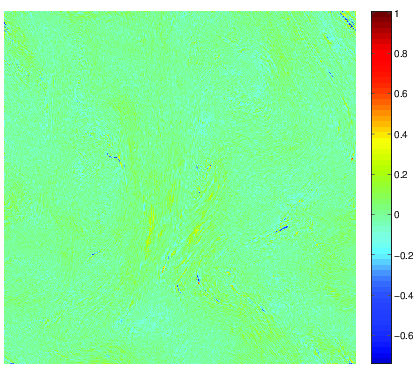}\\
(a)&(b)
\end{tabular}
\end{center}
\vspace{4mm}
\caption{(a): reconstructed data; (b): residual.}   \label{fig:plasma_rec}
\end{figure}

Figure \ref{fig:plasma_rec}(a) shows the image reconstructed by $L_1$ extension and Figure \ref{fig:plasma_rec}(b) 
gives the difference between the original image, Figure \ref{fig:plasma}(a) and the reconstructed image Figure \ref{fig:plasma_rec}(a).
As we can see, $L_1$ extension gives very good reconstruction. This result shows that the point integral method solve the isotropic elliptic equation 
very well on point cloud.

\section{Conclusion}

In this paper, we generalize the point integral method to solve the isotropic elliptic equation.
The point integral method is very easy to implement on point cloud, since it only needs the point cloud without any 
extra information. Moreover, it also has very good theoretical property. The coercivity of the original elliptic operator is partially
preserved in the point integral method. Based on this property, the convergence is proved. 
 

One important implication is the spectral convergence of the point integral method on random samples. 
Suppose the points are obtained by sampling a manifold according to some probability distribution $p(\bx)$. 
In the point integral method, the eigenvalue problem
\begin{align}
\label{eq:eigen}
  \left\{\begin{array}
    {cc}
    -\frac{1}{p^2(\bx)}\diver(p^2(\bx)\nabla u(\bx))=\lambda u(\bx), & \bx\in \M,\\
    \frac{\p u}{\p \bn}(\bx)=0,&\bx\in \p\M,    
  \end{array}\right.
\end{align}
is discretized as
\begin{equation}
\invt\sum_{\bx_j \in P} R_t(\bx_i,\bx_j)(u_i - u_j)  = 
\lambda\sum_{\bx_j \in P} \bar{R}_t(\bx_i,\bx_j) u_j.
\label{eqn:dis-eigen}
\end{equation}
This discretization is closely related with the normalized graph laplacian. 
Based on the theoretical results in this paper, it can be proved that the spectra of \eqref{eqn:dis-eigen}
converges to the spectra of \eqref{eq:eigen} as the number of sample points goes to infinity. 

The other interesting problem is how to generalize the point integral method to anisotropic elliptic equation. On this problem,
we already get some results. They are going to be reported in the subsequent paper.

\vspace{0.5in}
\noindent

\appendix

\bibliographystyle{abbrv}
\bibliography{reference}


\end{document}

%% file: integral_eqn.tex
Let $r(\bx)=L_t u-L_tu_t$ where $u$ and $u_t$ are the solution of \eqref{eq:neumann-homo} and \eqref{eq:integral-homo} respectively. 
Using integration by parts, we have
\begin{align}
\label{eq:r1}
  r(\bx)=&\frac{1}{t}\int_{\M}R_t(\bx,\by)(u(\bx)-u(\by))p(\by)\mathd\mu_\by-\int_{\M}\diver(p^2(\by)\nabla u(\by))\frac{\bar{R}_t(\bx,\by)}{p(\by)}\mathd\mu_\by\\
& \hspace{0cm}-2\int_{\p\M}\bar{R}_t(\bx,\by)\frac{\p u}{\p \bn}(\by)p(\by)\mathd\tau_\by\nonumber\\
=&\frac{1}{t}\int_{\M}(u(\bx)-u(\by)-(\bx-\by)\cdot\nabla u(\by))R_t(\bx,\by)p(\by)\mathd\mu_\by\nonumber\\
&\hspace{0cm}
 -\int_\M\Delta_\M u(\by)\bar{R}_t(\bx,\by)p(\by)\mathd\mu_\by\nonumber
\end{align}
The main idea of the proof is the Taylor expansion,
\begin{equation*}
  u(\bx)-u(\by)-(\bx-\by)\cdot\nabla u(\by)=\frac{1}{2}(\bx-\by)^T\cdot\mathbf{H}_u(\by)\cdot (\bx-\by)+O(|\bx-\by|^3)
\end{equation*}
where $\mathbf{H}_u(\by)$ is the Hessian matrix of $u$ at $\by$. 

Using the integration by parts, the second order term actually gives Laplace-Beltrami operator which cancel with the second term in \eqref{eq:r1}. 

In manifold, the Taylor expansion and integration by parts are more complicated. 
To make the whole idea rigorous, we need to introduce a special parametrization of the manifold $\M$. This parametrization is based on following proposition.
\begin{proposition}
Assume both $\M$ and $\p \M$ are $C^2$ smooth and
$\sigma$ is the minimum of the reaches of $\M$ and $\p \M$. 
For any point $\bx\in \M$, there is a neighborhood $U\subset \M$ of $\bx$, 
so that there is a parametrization 
$\Phi: \Omega \subset \R^k \rightarrow U$ satisfying the following conditions. For any $\rho \leq 0.1$, 
\begin{enumerate}
\item[(i)] $\Omega$ is convex and contains at least half of the ball
$B_{\Phi^{-1}(\bx)}(\frac{\rho}{5} \sigma)$, 
i.e., $vol(\Omega\cap B_{\Phi^{-1}(\bx)}(\frac{\rho}{5} \sigma)) > \frac{1}{2}(\frac{\rho}{5}\sigma)^k w_k$ where $w_k$ is the volume of unit ball in $\R^k$;
\item[(ii)] $B_{\bx}(\frac{\rho}{10} \sigma) \cap \M \subset U$. 
\item[(iii)] The determinant the Jacobian of $\Phi$ is bounded:
$(1-2\rho)^k \leq |D\Phi|  \leq (1+2\rho)^k$ over $\Omega$.  
\item[(iv)] For any points $\by, \bz \in U$, 
$1-2\rho \leq \frac{|\by-\bz|}{\left|\Phi^{-1}(\by) - \Phi^{-1}(\bz)\right|}  \leq 1+3\rho$.
\end{enumerate}
\label{prop:local_param}
\end{proposition}
This proposition basically says there exists a local parametrization of small distortion
if $(\M, \p\M)$ satisfies certain smoothness, and moreover, the parameter domain is 
convex and big enough. The proof of this proposition can be found in \cite{SS-neumann} and for the sake of completeness, 
we give the proof in the supplementary material. Next, we introduce a special parametrization of the manifold $\M$.

Let $\rho=0.1$, $\sigma$ be the minimum of the reaches of $\M$ and $\p \M$ and $\delta=\rho \sigma/20$. For any $\bx\in \M$, denote
\begin{eqnarray}
\label{eq:net}
  B_{\bx}^\delta=\left\{\by\in \mathcal{M}: |\bx-\by|\le \delta\right\},\quad \mathcal{M}_{\bx}^t=\left\{\by\in \mathcal{M}: |\bx-\by|^2\le 4t\right\}
\end{eqnarray}
and we assume $t$ is small enough such that $2\sqrt{t}\le \delta$.

Since the manifold $\mathcal{M}$ is compact, there exists a $\delta$-net, $\mathcal{N}_\delta=\{ \bq_i\in \mathcal{M},\;i=1,\cdots,N\}$, such that
\begin{eqnarray}
  \mathcal{M}\subset \bigcup_{i=1}^N B_{\bq_i}^\delta.\nonumber
\end{eqnarray}
and there exists a partition of $\M$, $\{\mathcal{O}_i, \;i=1,\cdots,N\}$, such that $\mathcal{O}_i\cap \mathcal{O}_j=\emptyset,\; i\neq j$ and
\begin{equation*}
  \M=\bigcup_{i=1}^N \mathcal{O}_i,\quad \mathcal{O}_i\subset B_{\bq_i}^\delta,\quad i=1,\cdots,N.
\end{equation*}

Using Proposition \ref{prop:local_param}, there exist
a parametrization $\Phi_i: \Omega_i\subset\mathbb{R}^k \rightarrow U_i\subset \mathcal{M},\; i=1,\cdots, N$, such that
\begin{itemize}
\item[1.] (Convexity) $B_{\bq_i}^{2\delta}\subset U_i$ and $\Omega_i$ is convex.
\item[2.] (Smoothness) $\Phi_i\in C^3(\Omega_i)$;
\item[3.] (Locally small deformation) 
For any points $\theta_1, \theta_2\in \Omega_i$, 
$$\frac{1}{2}\left|\theta_1-\theta_2\right| \leq \left\|\Phi_i(\theta_1)-\Phi_i(\theta_2)\right\|  
\leq 2\left|\theta_1-\theta_2\right|.$$
\end{itemize}
Using the partition, $\{\mathcal{O}_i, \;i=1,\cdots,N\}$, for any $\by\in \M$, there exists unique $J(\by)\in \{1,\cdots,N\}$, 
such that 
\begin{equation}
\label{def:index-j}
\by\in \mathcal{O}_{J(\by)}\subset B_{\bq_{J(\by)}}^\delta.
\end{equation}
Moerover, using the condition, $2\sqrt{t}\le \delta$, we have $\mathcal{M}_{\by}^t \subset B_{\bq_{J(\by)}}^{2\delta}\subset U_{J(\by)}$.
Then $\Phi_{J(\by)}^{-1}(\bx)$ and $\Phi_{J(\by)}^{-1}(\by)$ are both well defined for any $\bx\in \mathcal{M}_{\by}^t$.

Now, we define an auxiliary function, $\eta(\bx,\by)$ for any $\by\in \M,\;\bx\in \mathcal{M}_{\by}^t$. 
Let 
\begin{equation}
\label{def:eta}
\xi(\bx,\by)=\Phi_{J(\by)}^{-1}(\bx)-\Phi_{J(\by)}^{-1}(\by)\in \mathbb{R}^k,\quad \eta(\bx,\by)=\xi(\bx,\by)\cdot \p \Phi_{J(\by)}(\alpha(\bx,\by))\in \mathbb{R}^d, 
\end{equation}
where $\alpha(\bx,\by)=\Phi^{-1}_{J(\by)}(\by)$ and 
$\p$ is the gradient operator in the parameter space, i.e.
\begin{equation*}
  \p \Phi_j(\theta)=\left(\frac{\p \Phi_j}{\p \theta_1}(\theta),\frac{\p \Phi_j}{\p \theta_2}(\theta),\cdots,\frac{\p \Phi_j}{\p \theta_k}(\theta)\right),
\quad \theta\in \Omega_j\subset \mathbb{R}^k.
\end{equation*}


Now we state the proof of Theorem \ref{thm:integral_error}.
\begin{proof}

First, we split the residual $r(\bx)$ in \eqref{eq:r1} to four terms 
\begin{align}
  r(\bx)
=&r_1(\bx)+r_2(\bx)+r_3(\bx)-r_4(\bx)\nonumber
\end{align}
where
\begin{eqnarray}
r_1(\bx)&=&\frac{1}{t}\int_{\M}\left(u(\bx)-u(\by)-(\bx-\by)\cdot \nabla u(\by)-\frac{1}{2}\eta^i\eta^j(\nabla^i\nabla^j u(\by))\right)R_t(\bx,\by)p(\by)\mathd\mu_\by,\nonumber\\
r_2(\bx)&=&\frac{1}{2t}\int_{\M}\eta^i\eta^j(\nabla^i\nabla^j u(\by))R_t(\bx,\by)p(\by)\mathd\mu_\by
-\int_\M \eta^i(\nabla^i\nabla^j u(\by)\nabla^j\bar{R}_t(\bx,\by)p(\by)\mathd \mu_\by,\nonumber\\
r_3(\bx)&=& \int_\M \eta^i(\nabla^i\nabla^j u(\by)\nabla^j\bar{R}_t(\bx,\by)p(\by)\mathd \mu_\by
+\int_\M \mbox{div} \; \left(\eta^i(\nabla^i\nabla u(\by)\right)
\bar{R}_t(\bx,\by)p(\by)\mathd \mu_\by,\nonumber\\
r_4(\bx)&=&\int_\M \mbox{div} \; \left(\eta^i(\nabla^i\nabla u(\by)\right)\bar{R}_t(\bx,\by)p(\by)\mathd \mu_\by
+ \int_{\M}\Delta_\M u(\by)\bar{R}_t(\bx,\by)p(\by)\mathd \mu_\by.\nonumber
\nonumber
\end{eqnarray}
where $\nabla^i,\; i=1,\cdots,d$ is the $i$th component of the gradient $\nabla$, $\eta^i, \; i=1,\cdots,d$ is the $i$th component of $\eta(\bx,\by)$ defined in \eqref{def:eta}. 
To simplify the notation, 
we drop the variable $(\bx,\by)$ in the function $\eta(\bx,\by)$.

Next, we will prove the theorem by estimating above four terms one by one. 
First, we consider $r_1$. Let
$$d(\bx,\by)=u(\bx)-u(\by)-(\bx-\by)\cdot \nabla u(\by)-\frac{1}{2}\eta^i\eta^j(\nabla^i\nabla^j u(\by)).$$
we have
\begin{eqnarray}
  \int_\M |r_1(\bx)|^2\mathd\mu_\bx&=&  \int_\M \left|\int_\M R_t(\bx,\by)d(\bx,\by)p(\by)\mathd \mu_\by\right|^2\mathd\mu_\bx\nonumber\\
&\le & (\max_\by p(\by))^2 \int_\M \left(\int_\M R_t(\bx,\by)\mathd \mu_\by\right)\left(\int_\M R_t(\bx,\by)|d(\bx,\by)|^2\mathd \mu_\by\right)\mathd\mu_\bx\nonumber\\
&\le & C\int_\M \int_\M R_t(\bx,\by)|d(\bx,\by)|^2\mathd \mu_\by\mathd\mu_\bx\nonumber
\end{eqnarray}
and
\begin{eqnarray}
 \int_\M \int_\M R_t(\bx,\by)|d(\bx,\by)|^2\mathd \mu_\by\mathd\mu_\bx&=&
 \sum_{i=1}^N \int_\M\int_{\mathcal{O}_i}R_t(\bx,\by)|d(\bx,\by)|^2\mathd \mu_\by\mathd\mu_\bx\nonumber\\
&=&\sum_{i=1}^N \int_{\mathcal{O}_i}\left(\int_{\M_{\by}^t}R_t(\bx,\by)|d(\bx,\by)|^2\mathd \mu_\bx\right)\mathd\mu_\by.\nonumber
\end{eqnarray}
Using Newton-Leibniz formula, we get
\begin{eqnarray}
  d(\bx,\by)&=&u(\bx)-u(\by)-(\bx-\by)\cdot \nabla u(\by)-\frac{1}{2}\eta^i\eta^j(\nabla^i\nabla^j u(\by))\nonumber\\
&=&\xi^{i}\xi^{i'}\int_0^1\int_0^1\int_0^1s_1\frac{d}{d s_3}\left(\p_{i}\Phi^j(\al+s_3 s_1\xi)\p_{i'}\Phi^{j'}(\al+s_3s_2 s_1\xi)\nabla^{j'}\nabla^ju(\Phi(\al+s_3s_2s_1 \xi))\right)
\mathd s_3\mathd s_2\mathd s_1\nonumber\\
&=&\xi^{i}\xi^{i'}\xi^{i''}\int_0^1\int_0^1\int_0^1s_1^2s_2\p_{i}\Phi^j(\al+s_3 s_1\xi)\p_{i''}\p_{i'}\Phi^{j'}(\al+s_3s_2 s_1\xi)\nabla^{j'}\nabla^ju(\Phi(\al+s_3s_2s_1 \xi))
\mathd s_3\mathd s_2\mathd s_1\nonumber\\
&&+\xi^{i}\xi^{i'}\xi^{i''}\int_0^1\int_0^1\int_0^1s_1^2\p_{i''}\p_{i}\Phi^j(\al+s_3 s_1\xi)\p_{i'}\Phi^{j'}(\al+s_3s_2 s_1\xi)\nabla^{j'}\nabla^ju(\Phi(\al+s_3s_2s_1 \xi))
\mathd s_3\mathd s_2\mathd s_1\nonumber\\
&&+\xi^{i}\xi^{i'}\xi^{i''}\int_0^1\int_0^1\int_0^1s_1^2s_2\p_{i}\Phi^j(\al+s_3s_2 s_1\xi)\p_{i'}\Phi^{j'}(\al+s_3s_2 s_1\xi)\p_{i''}\Phi^{j''}(\al+s_3s_2 s_1\xi)\nonumber\\
&&\hspace{4cm}\nabla^{j''}\nabla^{j'}\nabla^ju(\Phi(\al+s_3s_2s_1 \xi))\mathd s_3\mathd s_2\mathd s_1\nonumber
\end{eqnarray}
Here, $\Phi^i,\; i=1,\cdots, d$ is the $i$th component of the parameterization function $\Phi$ and the parameterization function $\Phi=\Phi_{J(\by)}$, $J(\by)$ is the index function 
given in \eqref{def:index-j}. $\alpha=\alpha(\bx,\by)=\Phi_{J(\by)}^{-1}(\by)$, $\xi=\xi(\bx,\by)=\Phi_{J(\by)}^{-1}(\bx)-\Phi_{J(\by)}^{-1}(\by)$.
In the rest of the proof, without introducing any confusion, we always to use these short notations to save the space. 
In above derivation, we need the convexity property of the parameterization function to make sure all the integrals are well defined.

Using above equality and the smoothness of the parameterization functions, it is easy to show that
\begin{eqnarray}
&&  \int_{\mathcal{O}_i}\left(\int_{\M_{\by}^t}R_t(\bx,\by)|d(\bx,\by)|^2\mathd \mu_\bx\right)\mathd\mu_\by\nonumber\\
&\le & Ct^3 \int_0^1\int_0^1\int_0^1\int_{\mathcal{O}_i}\int_{\M_{\by}^t}R_t(\bx,\by)\left|D^{2,3}u(\Phi_{J(\by)}(\al+s_3s_2s_1 \xi))\right|^2\mathd \mu_\bx\mathd\mu_\by
\mathd s_3\mathd s_2\mathd s_1\nonumber\\
&\le & Ct^3 \max_{0\le s\le 1}\int_{\mathcal{O}_i}\int_{\M_{\by}^t}R_t(\bx,\by)\left|D^{2,3}u(\Phi_{i}(\al+s \xi))\right|^2\mathd \mu_\bx\mathd\mu_\by,\nonumber
\end{eqnarray}
where we use the fact that $J(\by)=i,\; \by \in \mathcal{O}_i$ and 
\begin{eqnarray*}
  \left|D^{2,3}u(\bx)\right|^2=\sum_{j,j',j''=1}^d|\nabla^{j''}\nabla^{j'}\nabla^ju(\bx)|^2
+\sum_{j,j'=1}^d|\nabla^{j'}\nabla^ju(\bx)|^2.
\end{eqnarray*}

Let $\bz_i=\Phi_i(\al+s \xi),\; 0\le s\le 1$, then for any $\by\in \mathcal{O}_i\subset B_{\bq_i}^\delta$ and $\bx\in \M_{\by}^t$,
\begin{eqnarray*}
  |\bz_i-\by|\le 2s|\xi|\le 4s|\bx-\by|\le 8s\sqrt{t},\quad |\bz_i-\bq_i|\le |\bz_i-\by|+|\by-\bq_i|\le \delta+8s\sqrt{t}.
\end{eqnarray*}
We can assume that $t$ is small enough such that $8\sqrt{t}\le \delta$, then we have
\begin{eqnarray*}
  \bz_i\in B_{\bq_i}^{2\delta}.
\end{eqnarray*}
After changing of variable, we obtain
\begin{eqnarray}
&&  \int_{\mathcal{O}_i}\int_{\M_{\by}^t}R_t(\bx,\by)\left|D^{2,3}u(\Phi_i(\al+s \xi))\right|^2\mathd \mu_\bx\mathd\mu_\by\nonumber\\
&\le & \frac{C}{\delta_0} \int_{\mathcal{O}_i}\int_{B_{\bq_i}^{2\delta}}\frac{1}{s^k}R\left(\frac{|\bz_i-\by|^2}{128s^2t}\right)
\left|D^{2,3}u(\bz_i)\right|^2\mathd \mu_{\bz_i}\mathd\mu_\by\nonumber\\
&=&\frac{C}{\delta_0} \int_{\mathcal{O}_i}\frac{1}{s^k}R\left(\frac{|\bz_i-\by|^2}{128s^2t}\right)\mathd\mu_\by
\int_{B_{\bq_i}^{2\delta}}\left|D^{2,3}u(\bz_i)\right|^2\mathd \mu_{\bz_i}\nonumber\\
&\le & C \int_{B_{\bq_i}^{2\delta}}\left|D^{2,3}u(\bx)\right|^2\mathd \mu_{\bx}.\nonumber
\end{eqnarray}
This estimate would give us that
\begin{eqnarray}
\label{est:r1}
  \|r_1(\bx)\|_{L^2(\M)}\le Ct^{1/2}\|u\|_{H^3(\M)}
\end{eqnarray}
Now, we turn to estimate the gradient of $r_1$.
\begin{eqnarray*}
  \int_\M |\nabla_\bx r_1(\bx)|^2\mathd\mu_\bx&\le &  C\int_\M \left|\int_\M \nabla_\bx R_t(\bx,\by)d(\bx,\by)p(\by)\mathd \mu_\by\right|^2\mathd\mu_\bx\\
&&+C\int_\M \left|\int_\M  R_t(\bx,\by)\nabla_\bx d(\bx,\by)p(\by)\mathd \mu_\by\right|^2\mathd\mu_\bx.
\end{eqnarray*}
where $\nabla_\bx$ is the gradient in $\M$ with respect to $\bx$.

Using the same techniques in the calculation of $\|r_1(\bx)\|_{L^2(\M)}$, we get that
the first term of right hand side can bounded as follows
\begin{eqnarray}
  \int_\M \left|\int_\M \nabla_\bx R_t(\bx,\by)d(\bx,\by)p(\by)\mathd \mu_\by\right|^2\mathd\mu_\bx\le C  \|u\|_{H^3(\M)}^2.\nonumber
\end{eqnarray}
The estimation of second term is a little involved. First, we have
\begin{eqnarray}
\int_\M \left|\int_\M  R_t(\bx,\by)\nabla_\bx d(\bx,\by)p(\by)\mathd \mu_\by\right|^2\mathd\mu_\bx
&\le & C\int_\M \left(\int_\M R_t(\bx,\by)|\nabla_\bx d(\bx,\by)|^2\mathd \mu_\by\right)\mathd\mu_\bx\nonumber\\
&=& C\sum_{i=1}^N \int_{\mathcal{O}_i}\left(\int_{\M_{\by}^t}R_t(\bx,\by)|\nabla_\bx d(\bx,\by)|^2\mathd \mu_\bx\right)\mathd\mu_\by.\nonumber
\end{eqnarray}
Also using Newton-Leibniz formula, we have
\begin{eqnarray}
  d(\bx,\by)
&=&\xi^{i}\xi^{i'}\int_0^1\int_0^1s_1\left(\p_{i}\Phi^j(\al+ s_1\xi)\p_{i'}\Phi^{j'}(\al+s_2 s_1\xi)\nabla^{j'}\nabla^ju(\Phi(\al+s_2s_1 \xi))\right)
\mathd s_2\mathd s_1\nonumber\\
&&-\xi^{i}\xi^{i'}\int_0^1\int_0^1s_1\left(\p_{i}\Phi^j(\al)\p_{i'}\Phi^{j'}(\al)\nabla^{j'}\nabla^ju(\Phi(\al))\right)
\mathd s_2\mathd s_1\nonumber
\end{eqnarray}
Then the gradient of $d(\bx,\by)$ has following representation,
\begin{eqnarray}
\nabla_\bx  d(\bx,\by)
&=&\xi^{i}\xi^{i'}\nabla_\bx\left(\int_0^1\int_0^1s_1\left(\p_{i}\Phi^j(\al+ s_1\xi)\p_{i'}\Phi^{j'}(\al+s_2 s_1\xi)\nabla^{j'}\nabla^ju(\Phi(\al+s_2s_1 \xi))\right)
\mathd s_2\mathd s_1\right)\nonumber\\
&&\hspace{-20mm}+\nabla_\bx\left(\xi^{i}\xi^{i'}\right)\int_0^1\int_0^1\int_0^1s_1\frac{d}{d s_3}\left(\p_{i}\Phi^j(\al+ s_3s_1\xi)\p_{i'}\Phi^{j'}(\al+s_3s_2 s_1\xi)
\nabla^{j'}\nabla^ju(\Phi(\al+s_3s_2s_1 \xi))\right)
\mathd s_3\mathd s_2\mathd s_1\nonumber\\
&=&d_1(\bx,\by)+d_2(\bx,\by).\nonumber
\end{eqnarray}
For $d_1$, we have
\begin{eqnarray}
  \int_{\mathcal{O}_i}\left(\int_{\M_{\by}^t}R_t(\bx,\by)|d_1(\bx,\by)|^2\mathd \mu_\bx\right)\mathd\mu_\by
&\le& Ct^2\max_{0\le s\le 1} \int_{\mathcal{O}_i}\left(\int_{\M_{\by}^t}R_t(\bx,\by)|D^{2,3}u(\Phi_i(\al+s \xi))|^2\mathd \mu_\bx\right)\mathd\mu_\by,\nonumber
\end{eqnarray}
which means that
\begin{eqnarray}
\label{eqn:est-dr1-d1}
   \int_{\mathcal{O}_i}\left(\int_{\M_{\by}^t}R_t(\bx,\by)|d_1(\bx,\by)|^2\mathd \mu_\bx\right)\mathd\mu_\by\le C  \int_{B_{\bq_i}^{2\delta}}|D^{2,3}u(\bx)|^2\mathd \mu_\bx
\end{eqnarray}
For $d_2$, we have
\begin{eqnarray}
&&  d_2(\bx,\by)\nonumber\\
  &=&\nabla_\bx\left(\xi^{i}\xi^{i'}\right)\int_{[0,1]^3}s_1\frac{d}{d s_3}\left(\p_{i}\Phi^j(\al+ s_3s_1\xi)\p_{i'}\Phi^{j'}(\al+s_3s_2 s_1\xi)
\nabla^{j'}\nabla^ju(\Phi(\al+s_3s_2s_1 \xi))\right)
\mathd s_3\mathd s_2\mathd s_1\nonumber\\
&=&\nabla_\bx\left(\xi^{i}\xi^{i'}\right)\xi^{i''}\int_{[0,1]^3}s_1^2s_2\p_{i}\Phi^j(\al+s_3 s_1\xi)\p_{i''}\p_{i'}\Phi^{j'}(\al+s_3s_2 s_1\xi)\nabla^{j'}\nabla^ju(\Phi(\al+s_3s_2s_1 \xi))
\mathd s_3\mathd s_2\mathd s_1\nonumber\\
&&+\nabla_\bx\left(\xi^{i}\xi^{i'}\right)\xi^{i''}\int_{[0,1]^3}s_1^2\p_{i''}\p_{i}\Phi^j(\al+s_3 s_1\xi)\p_{i'}\Phi^{j'}(\al+s_3s_2 s_1\xi)\nabla^{j'}\nabla^ju(\Phi(\al+s_3s_2s_1 \xi))
\mathd s_3\mathd s_2\mathd s_1\nonumber\\
&&+\nabla_\bx\left(\xi^{i}\xi^{i'}\right)\xi^{i''}\int_{[0,1]^3}s_1^2s_2\p_{i}\Phi^j(\al+s_2 s_1\xi)\p_{i'}\Phi^{j'}(\al+s_3s_2 s_1\xi)\p_{i''}\Phi^{j''}(\al+s_3s_2 s_1\xi)\nonumber\\
&&\hspace{4cm}\nabla^{j''}\nabla^{j'}\nabla^ju(\Phi(\al+s_3s_2s_1 \xi))\mathd s_3\mathd s_2\mathd s_1\nonumber
\end{eqnarray}
This formula tells us that
\begin{eqnarray}
 \int_{\mathcal{O}_i}\left(\int_{\M_{\by}^t}R_t(\bx,\by)|d_2(\bx,\by)|^2\mathd \mu_\bx\right)\mathd\mu_\by
&\le& Ct^2\max_{0\le s\le 1} \int_{\mathcal{O}_i}\left(\int_{\M_{\by}^t}R_t(\bx,\by)|D^{2,3}u(\Phi(\al+s \xi))|^2\mathd \mu_\bx\right)\mathd\mu_\by.\nonumber
\end{eqnarray}
Using the same arguments as that in the calculation of $\|r_1\|_{L^2(\M)}$, we have
\begin{eqnarray}
\label{eqn:est-dr1-d2}
   \int_{\mathcal{O}_i}\left(\int_{\M_{\by}^t}R_t(\bx,\by)|d_2(\bx,\by)|^2\mathd \mu_\bx\right)\mathd\mu_\by\le C  \int_{B_{\bq_i}^{2\delta}}|D^3u(\bx)|^2\mathd \mu_\bx
\end{eqnarray}
Combining \eqref{eqn:est-dr1-d1} and \eqref{eqn:est-dr1-d2}, we have
\begin{eqnarray}
\label{est:dr1}
  \|\nabla r_1(\bx)\|_{L^2(\M)}\le C\|u\|_{H^3(\M)}
\end{eqnarray}

For $r_2$, first, notice that
\begin{eqnarray}
  \nabla^j\bar{R}_t(\bx,\by)&=&\frac{1}{2t}\p_{m'}\Phi^j(\al) g^{m'n'}\p_{n'}\Phi^i(\al) (x^i-y^i)\hk,\nonumber\\
\frac{\eta^j}{2t}R_t(\bx,\by)&=&\frac{1}{2t}\p_{m'}\Phi^j(\al) g^{m'n'}\p_{n'}\Phi^i(\al) \xi^{i'}\p_{i'}\Phi^i\hk.\nonumber
\end{eqnarray}
Then, we have
\begin{eqnarray}
  &&\nabla^j\bar{R}_t(\bx,\by)-\frac{\eta^j}{2t}R_t(\bx,\by)\nonumber\\
&=&\frac{1}{2t}\p_{m'}\Phi^i g^{m'n'}\p_{n'}\Phi^j \left(x^j-y^j-\xi^{i'}\p_{i'}\Phi^j\right)\hk\nonumber\\
&=&\frac{1}{2t}\xi^{i'}\xi^{j'}\p_{m'}\Phi^i g^{m'n'}\p_{n'}\Phi^j \left(\int_0^1\int_0^1s\p_{j'}\p_{i'}\Phi^j(\al+\tau s \xi)\mathd \tau\mathd s\right)\hk\nonumber
\end{eqnarray}
Thus, we get
\begin{eqnarray}
 \left|\nabla^j\bar{R}_t(\bx,\by)-\frac{\eta^j}{2t}R_t(\bx,\by)\right|&\le& \frac{C|\xi|^2}{t}\hk\nonumber\\
\left|\nabla_\bx\left(\nabla^j\bar{R}_t(\bx,\by)-\frac{\eta^j}{2t}R_t(\bx,\by)\right)\right|&\le& \frac{C|\xi|}{t}\hk+\frac{C|\xi|^3}{t^2}|R'_t(\bx,\by)|\nonumber
\end{eqnarray}
Then, we have following bound for $r_2$,
\begin{align}
\label{est:r2}
 & \int_{\M}|r_2(\bx)|^2\mathd \mu_\bx\\
\le& Ct\int_\M  \left(\int_\M \hk |D^2u(\by)|p(\by)\mathd \mu_\by\right)^2\mathd\mu_\bx\nonumber\\
\le &  Ct\int_\M  \left(\int_\M \hk p(\by) \mathd\mu_\by\right)\int_\M \hk |D^2u(\by)|^2p(\by)\mathd \mu_\by\mathd\mu_\bx \nonumber\\
\le& Ct  \max_{\by}\left(\int_\M \hk\mathd \mu_\bx\right) \int_\M|D^2u(\by)|^2p(\by)\mathd \mu_\by \nonumber\\
\le& Ct\|u\|_{H^2(\M)}^2.\nonumber
\end{align}
Similarly, we have
\begin{align}
\label{est:dr2}
 &\int_M |\nabla r_2(\bx)|^2\mathd\mu_\bx\\
\le &  Ct\int_\M  \left(\int_\M\nabla_\bx \hk p(\by)\mathd\mu_\by\right)\int_\M \nabla_\bx\hk
|D^2u(\by)|^2p(\by)\mathd \mu_\by\mathd\mu_\bx\nonumber\\
 \le& C\sqrt{t} 
\max_{\by}\left(\int_\M \nabla_\bx\hk\mathd \mu_\bx\right) \int_\M|D^2u(\by)|^2p(\by)\mathd \mu_\by\nonumber\\
\le&  C\|u\|_{H^2(\M)}^2. \nonumber
\end{align}
$r_3$ is relatively easy to estimate by using the well known Gauss formula.
\begin{align*}
  r_3(\bx)&=\int_{\p\mathcal{M}}n^{j}\eta^{i}(\nabla^{i} \nabla^{j}u(\by)) \rhk p(\by) \mathd\tau_\by
-\int_{\mathcal{M}}\eta^{i}(\nabla^{i} \nabla^{j}u(\by)) \rhk \nabla^{j}p(\by)\mathd\mu_\by\\
&=\tilde{I}_{bd}-\int_{\mathcal{M}}\eta^{i}(\nabla^{i} \nabla^{j}u(\by)) \rhk \nabla^{j}p(\by)\mathd\mu_\by
\end{align*}
where $\tilde{I}_{bd}=\int_{\p\mathcal{M}}n^{j}\eta^{i}(\nabla^{i} \nabla^{j}u(\by)) \rhk p(\by) \mathd\tau_\by$.

Using the assumption that $p\in C^1(\M)$, it is easy to get that 
\begin{align}
  \label{est:r3}
  \|r_3-\tilde{I}_{bd}\|_{L^2(\M)}&\le C\sqrt{t}\|u\|_{H^2(\M)},\\
\label{est:dr3}
\|\nabla(r_3-\tilde{I}_{bd})\|_{L^2(\M)}&\le C\|u\|_{H^2(\M)}.
\end{align}

Now, we turn to bound the last term $r_4$. Notice that
\begin{eqnarray}
\label{eqn:grad2Delta} 
 \nabla^j\left(\nabla^{j}u(\by)\right)&=&(\p_{k'} \Phi^j)g^{k'l'}\p_{l'}\left((\p_{m'}\Phi^j)g^{m'n'}(\p_{n'} u )\right)\\
&=&(\p_{k'}\Phi^j)g^{k'l'}\left(\p_{l'}(\p_{m'}\Phi^j)\right)g^{m'n'}(\p_{n'}u)\nonumber\\
&&+(\p_{k'}\Phi^j)g^{k'l'}(\p_{m'}\Phi^j)\p_{l'}\left(g^{m'n'}(\p_{n'}u)\right)\nonumber\\
&=&\frac{1}{\sqrt{\det G}}(\p_{m'}\sqrt{\det G}) g^{m'n'}(\p_{n'}u)+\p_{m'}\left(g^{m'n'}(\p_{n'}u)\right)\nonumber\\
&=&\frac{1}{\sqrt{\det G}}\p_{m'}\left(\sqrt{\det G} g^{m'n'}(\p_{n'}u)\right) =\Delta_\M u(\by).\nonumber
\end{eqnarray}
where $\det G$ is the determinant of $G$ and $G=(g_{ij})_{i,j=1,\cdots,k}$.
Here we use the fact that
\begin{eqnarray}
  (\p_{k'}\Phi^j)g^{k'l'}\left(\p_{l'}(\p_{m'}\Phi^j)\right)&=&(\p_{k'}\Phi^j)g^{k'l'}\left(\p_{m'}(\p_{l'}\Phi^j)\right)\nonumber\\
&=&(\p_{m'}(\p_{k'}\Phi^j))g^{k'l'}(\p_{l'}\Phi^j)\nonumber\\
&=&\frac{1}{2}g^{k'l'}\p_{m'}(g_{k'l'})\nonumber \\
&=&\frac{1}{\sqrt{\det G}}(\p_{m'}\sqrt{\det G}).\nonumber
\end{eqnarray}
Moreover, we have
\begin{eqnarray}
\label{eqn:div}
&& g^{i'j'}(\p_{j'}\Phi^{j})(\p_{i'}\xi^{l})(\p_{l}\Phi^{i})(\nabla^{i}\nabla^{j}u(\by))\\
&=&- g^{i'j'}(\p_{j'}\Phi^{j})(\p_{i'}\Phi^{i})(\nabla^{i}\nabla^{j}u(\by))\nonumber\\
&=&- g^{i'j'}(\p_{j'}\Phi^{j})(\p_{i'}\Phi^{i})(\p_{m'}\Phi^{i})g^{m'n'}\p_{n'}\left(\nabla^{j}u(\by) \right)\nonumber\\
&=&- g^{i'j'}(\p_{j'}\Phi^{j})\p_{i'}\left(\nabla^{j}u(\by)\right)\nonumber\\
&=&- \nabla^{j}\left(\nabla^{j}u(\by)\right).\nonumber
\end{eqnarray}
where the first equalities are due to that
$\p_{i'}\xi^l = -\delta_{i'}^l$.
Then we have
\begin{eqnarray}
  &&\mbox{div} \; \left(\eta^i(\nabla^i\nabla^j u(\by))\right)+\Delta_\M u(\by)\nonumber\\
&=&\frac{1}{\sqrt{\det G}}\,\p_{i'}\left(\sqrt{\det G}\,g^{i'j'}(\p_{j'}\Phi^j) \xi^{l}(\p_l \Phi^i)(\nabla^i\nabla^j u(\by))\right)
-g^{i'j'}(\p_{j'}\Phi^{j})(\p_{i'}\xi^{l})(\p_{l}\Phi^{i})(\nabla^{i}\nabla^{j}u(\by))\nonumber\\
&=&\frac{\xi^l}{\sqrt{\det G}}\,\p_{i'}\left(\sqrt{\det G}\,g^{i'j'}(\p_{j'}\Phi^j) (\p_l \Phi^i)(\nabla^i\nabla^j u(\by))\right).\nonumber
\end{eqnarray}
Here we use the equalities \eqref{eqn:grad2Delta}, \eqref{eqn:div}, $\eta^i = \xi^{l} \p_{i'}\Phi^l$ and the definition of $\mbox{div}$,
\begin{equation}
\mbox{div}X = \frac{1}{\sqrt{\det G}}\p_{i'} (\sqrt{\det G}\,g^{i'j'}\p_{j'}\Phi^k X^k).
\label{eqn:divergence}
\end{equation}
where $X$ is a smooth tangent vector
field on $\M$ and $(X^1,\dots,X^d)^t$ is its representation in embedding coordinates.

Hence,
\begin{eqnarray}
  r_4(\bx)=\int_\mathcal{M}\frac{\xi^l}{\sqrt{\det G}}\,\p_{i'}\left(\sqrt{\det G}\,g^{i'j'}
(\p_{j'}\Phi^j) (\p_l \Phi^i)(\nabla^i\nabla^j u(\by))\right) \rhk p(\by)\mathd \mu_\by\nonumber
\end{eqnarray}
Then it is easy to get that
\begin{eqnarray}
  \label{est:r4}
  \|r_4(\bx)\|_{L^2(\M)}&\le& C t^{1/2}\|u\|_{H^3(\M)},\\
\|\nabla r_4(\bx)\|_{L^2(\M)}&\le & C\|u\|_{H^3(\M)}.
\label{est:dr4}
\end{eqnarray}
By combining \eqref{est:r1},\eqref{est:dr1},\eqref{est:r2},\eqref{est:dr2},\eqref{est:r3},\eqref{est:dr3},\eqref{est:r4},\eqref{est:dr4}, we know that
\begin{eqnarray}
  \label{est:r0}
  \|r-\tilde{I}_{bd}\|_{L^2(\M)}&\le& C t^{1/2}\|u\|_{H^3(\M)},\\
\|\nabla (r-\tilde{I}_{bd})\|_{L^2(\M)}&\le & C\|u\|_{H^3(\M)}.
\label{est:dr0}
\end{eqnarray}
Using the definition of $I_{bd}$ and $\tilde{I}_{bd}$, we obtain
\begin{equation*}
  I_{bd}-\tilde{I}_{bd}=\int_{\p\mathcal{M}}n^{j}(\by)(\bx-\by-\eta(\bx,\by))\cdot(\nabla \nabla^{j}u(\by)) \rhk p(\by) \mathd\tau_\by
\end{equation*}
Using the definition of $\eta(\bx,\by)$, it is easy to check that
\begin{equation*}
  |\bx-\by-\eta(\bx,\by)|=O(|\bx-\by|^2),\quad |\nabla_{\bx}(\bx-\by-\eta(\bx,\by))|=O(|\bx-\by|)
\end{equation*}
which implies that 
\begin{eqnarray}
  \label{est:bd}
  \|I_{bd}-\tilde{I}_{bd}\|_{L^2(\M)}&\le& C t^{3/4}\|u\|_{H^2(\M)},\\
\|\nabla (I_{bd}-\tilde{I}_{bd})\|_{L^2(\M)}&\le & Ct^{1/4}\|u\|_{H^3(\M)}.
\label{est:dbd}
\end{eqnarray}
The theorem is proved by putting \eqref{est:r0}, \eqref{est:dr0}, \eqref{est:bd}, \eqref{est:dbd} together.
\end{proof}


  

\begin{remark}
  Using above proof, we can also show that the $L_2$ error in the integral approximation \eqref{eq:integral-adaptive} is $O(t^{1/4})$. 
\end{remark}


%% file: discretization.tex
To simplify the notation, we introduce a intermediate operator defined as follows,
\begin{eqnarray}
 L_{t, h}u(\bx) &=& \invt\sum_{\bfp_j \in P} \hkxpj(u(\bx) - u(\bfp_j))p(\bx_j)V_j.
 \label{eqn:laplace_dis} 
\end{eqnarray}
Let  $u_{t,h}=I_{\mathbf{f}}(\mathbf{u})$ with $\mathbf{u}$
satisfying equation \eqref{eqn:dis-homo} and $I_{\bff}$ is given in \eqref{eqn:interp_neumann}. 
One can verify that the following equation are satisfied,
\begin{eqnarray}
-L_{t, h}u_{t,h}(\bx) &=& \sum_{\bfp_j \in P} \rhkxpj f(\bfp_j)/p(\bx_j)V_j.
\label{eqn:integral_dis_interp}
\end{eqnarray}

In the proof, we need {\it a prior} estimate of $\bfu$ which is given as following.
\begin{theorem}
Suppose $\mathbf{u}=(u_1, \cdots, u_{|P|})$ with $\sum_{i=1}^{|P|} u_i p_iV_i = 0$ solves the problem~\eqref{eqn:dis-homo}
and $\bff =(f(\bfp_1), \cdots, f(\bfp_{|P|}))^t$ for $f\in C(\mathcal{M})$. 
Then there exists a constant $C>0$ such that 
\begin{eqnarray*}
  \left(\sum_{i=1}^{|P|}u_i^2p_iV_i\right)^{1/2}\le C\|f\|_\infty,
\end{eqnarray*}
provided $t$ and $\frac{h(P,\mathbf{V},\M)}{\sqrt{t}}$ are small enough.
\label{lem:bound_solution_bfu}
\end{theorem}

This theorem is an easy corollary of following theorem.
\begin{theorem}
If the manifolds $\M$ is $C^\infty$, there exist constants $C>0,\, C_0>0$ independent on $t$ so that for 
any ${\bf u} = (u_1, \cdots, u_{|P|})^t \in \R^d$ with $\sum_{i=1}^{|P|} u_ip_iV_i = 0$ and for any sufficient small $t$ 
and $\frac{h(P,\mathbf{V},\M)}{\sqrt{t}}$ 
\begin{equation*}
\sum_{i,j=1}^{|P|}R_t(\bx_i,\bx_j)(u_i-u_j)^2p_ip_jV_iV_j  \geq C(1-\frac{C_0h(P,\mathbf{V},\M)}{\sqrt{t}}) \sum_{i=1}^{|P|} u_i^2p_iV_i.
\end{equation*}
\label{thm:elliptic_L}
\end{theorem}
The proof of this theorem is given in the supplementary material which is a small modification of the proof of Theorem 9.1 in \cite{SS-neumann}.



We are now ready to prove Theorem \ref{thm:dis_error}.
\begin{proof}
To simplify the notation, we denote $h=h(P,\mathbf{V},\M)$ and $n=|P|$
and denote
\begin{eqnarray}
u_{t,h}(\bx) = \frac{1}{w_{t,h}(\bx)}\left(\sum_{\bfp_j\in P}R_t(\bx,\bfp_j)u_jp_jV_j
-t\sum_{\bfp_j\in P}\bar{R}_t(\bx,\bfp_j)f_jV_j/p_j\right)
\end{eqnarray}
where $\mathbf{u}=(u_1, \cdots, u_n)^t$ with $\sum_{i=1}^n u_ip_i V_i = 0$ solves the problem~\eqref{eqn:dis-homo}, 
$ f_j = f(\bfp_j)$ and $w_{t,h}(\bx)=\sum_{\bfp_j\in P}R_t(\bx,\bfp_j)p_jV_j$.  
For convenience, we set
\begin{eqnarray}
  a_{t,h}(\bx)&=&\frac{1}{w_{t,h}(\bx)}\sum_{\bfp_j\in P}R_t(\bx,\bfp_j)u_jp_jV_j,\\
c_{t,h}(\bx)&=&-\frac{t}{w_{t,h}(\bx)}\sum_{\bfp_j\in P}\bar{R}_t(\bx,\bfp_j)f(\bfp_j)V_j/p_j, 
\end{eqnarray}
and thus $u_{t,h} = a_{t, h} +c_{t, h}$.

First we upper bound $\|L_t(u_{t,h}) - L_{t, h}(u_{t, h})\|_{L^2(\M)}$.
For $c_{t, h}$, we have 
\begin{eqnarray}
&&\left|\left(L_tc_{t, h} - L_{t, h}c_{t, h}\right)(\bx)\right|\nonumber \\
&=& \frac{1}{t} \left|\int_{\mathcal{M}}R_t(\bx,\by)(c_{t,h}(\bx)-c_{t,h}(\by))p(\by)  \mathd \mu_\by-\sum_{\bfp_j\in P}R_t(\bx,\bfp_j)(c_{t,h}(\bx)-c_{t,h}(\bfp_j))p_jV_j\right|\nonumber\\
&\le &\frac{1}{t} \left|c_{t,h}(\bx)\right|\left|\int_{\mathcal{M}}R_t(\bx,\by) p(\by)\mathd \mu_\by-\sum_{\bfp_j\in P}R_t(\bx,\bfp_j)p_jV_j\right|\nonumber\\
&&+ \frac{1}{t}\left|\int_{\mathcal{M}}R_t(\bx,\by)c_{t,h}(\by)  p(\by)\mathd \mu_\by-\sum_{\bfp_j\in P}R_t(\bx,\bfp_j)c_{t,h}(\bfp_j)p_jV_j\right|\nonumber\\
&\le &\frac{Ch}{t^{3/2}}\left|c_{t,h}(\bx)\right|+ \frac{Ch}{t^{3/2}} \|c_{t,h}\|_{C^1(\M)}\nonumber\\
&\le &\frac{Ch}{t^{3/2}}t\|f\|_\infty+ \frac{Ch}{t^{3/2}}(t\|f\|_\infty + t^{1/2}\|f\|_{\infty})\le\frac{Ch}{t}\|f\|_{\infty}.\nonumber
\end{eqnarray}
For $a_{t, h}$, we have 
\begin{eqnarray}
\label{eqn:a_1}
&&  \int_{\mathcal{M}}\left(a_{t,h}(\bx)\right)^2\left|\int_{\mathcal{M}}R_t(\bx,\by) p(\by)\mathd \mu_\by-\sum_{\bfp_j\in P}R_t(\bx,\bfp_j)p_jV_j\right|^2\mathd\mu_\bx\\
&\le & \frac{Ch^2}{t}\int_{\mathcal{M}}\left(a_{t,h}(\bx)\right)^2\mathd\mu_\bx 
\le  \frac{Ch^2}{t} \int_{\mathcal{M}} \left( \frac{1}{w_{t, h}(\bx)} \sum_{\bfp_j\in P}R_t(\bx,\bfp_j)u_jp_jV_j \right)^2 \mathd\mu_\bx \nonumber \\
&\le & \frac{Ch^2}{t} \int_{\mathcal{M}} \left( \sum_{\bfp_j\in P}R_t(\bx,\bfp_j)u_j^2p_jV_j  \right) \left( \sum_{\bfp_j\in P}R_t(\bx,\bfp_j)p_jV_j  \right) \mathd\mu_\bx \nonumber \\
&\le & \frac{Ch^2}{t} \left( \sum_{j=1}^{n}u_j^2p_jV_j \int_{\mathcal{M}}R_t(\bx,\bfp_j)  \mathd\mu_\bx    \right) \le \frac{Ch^2}{t}\sum_{j=1}^{n}u_j^2p_jV_j. \nonumber
\end{eqnarray}
Let 
\begin{eqnarray}
A &=&  C_t\int_{\mathcal{M}}\frac{1}{w_{t, h}(\by)}R\left(\frac{|\bx-\by|^2}{4t}\right)R\left(\frac{|\bfp_i-\by|^2}{4t}\right) p(\by)\mathd \mu_\by\nonumber\\
 &-&C_t\sum_{\bfp_j\in P} \frac{1}{w_{t, h}(\bfp_j)}R\left(\frac{|\bx-\bfp_j|^2}{4t}\right)R\left(\frac{|\bfp_i-\bfp_j|^2}{4t}\right)p_jV_j. \nonumber
\end{eqnarray}
We have $|A|<\frac{Ch}{t^{1/2}}$ for some constant $C$ independent of $t$. In addition, notice that
only when $|\bx-\bfp_i|^2\leq 16t $ is $A\neq 0$, which implies 
\begin{eqnarray}
|A| \leq \frac{1}{\delta_0}|A|R\left(\frac{|\bx-\bfp_i|^2}{32t}\right). \nonumber
\end{eqnarray}
Then we have
\begin{eqnarray}
\label{eqn:a_2}
&&\int_{\mathcal{M}}\left|\int_{\mathcal{M}}R_t(\bx,\by)a_{t,h}(\by)p(\by)  \mathd \mu_\by-\sum_{\bfp_j\in P}R_t(\bx,\bfp_j)a_{t,h}(\bfp_j)p_jV_j\right|^2\mathd\mu_\bx\\
&=& \int_{\mathcal{M}}\left(\sum_{i=1}^{n}C_tu_ip_iV_i A \right)^2\mathd\mu_\bx\le
\frac{Ch^2}{t} \int_{\mathcal{M}}\left(\sum_{i=1}^{n} C_t|u_i|p_iV_i R\left(\frac{|\bx-\bfp_i|^2}{32t}\right)  \right)^2 \mathd\mu_\bx \nonumber\\
&\le &\frac{Ch^2}{t} \int_{\mathcal{M}} \left(\sum_{i=1}^{n} C_t R\left(\frac{|\bx-\bfp_i|^2}{32t}\right) u^2_ip_iV_i\right) 
 \left(\sum_{\bfp_i\in P} C_t R\left(\frac{|\bx-\bfp_i|^2}{32t} \right)p_iV_i  \right) \mathd\mu_\bx \nonumber\\
&\le &\frac{Ch^2}{t} \sum_{i=1}^{n} \left(\int_{\mathcal{M}} C_t R\left(\frac{|\bx-\bfp_i|^2}{32t}\right)  \mathd\mu_\bx \left(u^2_ip_iV_i\right)\right)  \le
 \frac{Ch^2}{t} \left(\sum_{i=1}^{n}u_i^2p_iV_i\right). \nonumber
\end{eqnarray} 
Combining Equation~\eqref{eqn:a_1},~\eqref{eqn:a_2} and Lemma~\ref{lem:bound_solution_bfu}, 
\begin{eqnarray}
&&\|L_ta_{t, h} - L_{t, h}a_{t, h}\|_{L^2(\M)} \nonumber \\
&=&\left(\int_\M \left|\left(L_t(a_{t, h}) - L_{t, h}(a_{t, h})\right)(\bx)\right|^2 \mathd\mu_\bx\right)^{1/2} \nonumber\\
&\le & \frac{1}{t}\left(\int_{\mathcal{M}}\left(a_{t,h}(\bx)\right)^2\left|\int_{\mathcal{M}}R_t(\bx,\by)p(\by) \mathd \mu_\by
-\sum_{\bfp_j\in P}R_t(\bx,\bfp_j)p_jV_j\right|^2\mathd\mu_\bx\right)^{1/2}\nonumber\\
&&+  \frac{1}{t}\left(\int_{\mathcal{M}}\left|\int_{\mathcal{M}}R_t(\bx,\by)a_{t,h}(\by)  p(\by)\mathd \mu_\by-\sum_{\bfp_j\in P}R_t(\bx,\bfp_j)a_{t,h}(\bfp_j)p_jV_j\right|^2\mathd\mu_\bx
\right)^{1/2} \nonumber \\
&\le & \frac{Ch}{t^{3/2}} \left(\sum_{i=1}^{n}u_i^2p_iV_i\right)^{1/2} \leq  \frac{Ch}{t^{3/2}}\|f\|_\infty.\nonumber
\end{eqnarray}
Assembling the parts together, we have the following upper bound.
\begin{eqnarray}
\label{eqn:elliptic_L_1}
&&\|L_tu_{t, h} - L_{t, h}u_{t, h}\|_{L^2(\M)} \\
&\le& \|L_ta_{t, h} - L_{t, h}a_{t, h}\|_{L^2(\M)}  + \|L_tc_{t, h} - L_{t, h}c_{t, h}\|_{L^2(\M)} \nonumber \\
&\le& \frac{Ch}{t^{3/2}}\|f\|_\infty + \frac{Ch}{t}\|f\|_{\infty} \le \frac{Ch}{t^{3/2}} \|f\|_\infty\nonumber
\end{eqnarray}
At the same time, since $u_t$ respectively $u_{t,h}$ solves equation~\eqref{eq:integral-homo} respectively equation~\eqref{eqn:integral_dis_interp}, 
we have 
\begin{eqnarray}
\label{eqn:elliptic_L_2}
&&\|L_t(u_{t}) - L_{t, h}(u_{t, h})\|_{L^2(\M)}  \\ 
&=&\left(\int_\M \left(\left(L_tu_{t} - L_{t, h}u_{t, h}\right)(\bx)\right)^2     \mathd\mu_\bx\right)^{1/2}  \nonumber \\
&=& \left( \int _\M \left( \int_{\mathcal{M}}\bar{R}_t(\bx,\by)f(\by)/p(\by) - \sum_{\bfp_j\in P}\bar{R}_t(\bx,\bfp_j)f(\bfp_j)V_j/p_j\right)^2 \mathd\mu_\bx\right)^{1/2}\nonumber\\ 
&\le&  \frac{Ch}{t^{1/2}}\|f\|_{C^1(\M)}. \nonumber
\end{eqnarray}
The complete $L^2$ estimate follows from Equation~\eqref{eqn:elliptic_L_1} and~\eqref{eqn:elliptic_L_2}. 
 
The estimate of the gradient, $\|\nabla(L_t(u_{t}) - L_{t, h}(u_{t, h}))\|_{L^2(\M)}$, can be obtained similarly. 
\end{proof}
